%% file: marina2.tex
\documentclass [10pt,reqno]{amsart}
\usepackage {hannah}
\usepackage {color,graphicx}
\usepackage {multirow}
\usepackage {texdraw}
\usepackage [latin1]{inputenc}

\DeclareMathOperator {\ev}{ev}
\DeclareMathOperator {\ft}{f\/t}
\renewcommand{\Z}{\mathbb{Z}}

\renewcommand{\F}{\mathbb{F}}
\renewcommand{\P}{\mathbb{P}}
\newcommand{\M}{\mathcal{M}}
\newcommand{\numb}[4]{\mathcal{N}^{\emph{#1}}_{\mathbb{F}_{\emph{#2}}}(\emph{#3},\emph{#4})}
\newcommand{\numbc}[4]{N^{\emph{#1}}_{\mathbb{F}_{\emph{#2}}}(\emph{#3},\emph{#4})}

\title [Tropical enumerative invariants of $\F_0$ and $\F_2$]{Tropical enumerative invariants of $\F_0$ and $\F_2$}

\author {Marina Franz and Hannah Markwig}

\address{ Fachbereich Mathematik, TU Kaiserslautern, Postfach 3049, 67653 Kaiserslautern,
 Germany}
\email{franz@mathematik.uni-kl.de}
\address {CRC ``Higher Order Structures in Mathematics'', Georg-August-Universit\"at G\"ottingen, Bunsenstr. 3-5, D-37073 G\"ottingen, Germany}
\email {hannah@uni-math.gwdg.de}

\thanks{Supported by the German Research Foundation 
(Deutsche Forschungsgemeinschaft (DFG)) through 
the Institutional Strategy of the University of G\"ottingen}

\begin {document}

\begin {abstract}
 There is an equation relating numbers of curves on $\F_0=\P^1\times\P^1$
 satisfying incidence conditions and numbers of curves on $\F_2$ satisfying incidence conditions.
 The purpose of this paper is to give a tropical proof of this equation in the case of rational curves. We use induction on the degree and two Kontsevich-type formulas for curves on $\F_0$ and on $\F_2$. The formula for $\F_2$ was not known before and is proved using tropical geometry. 
\end {abstract}

\maketitle

\section{Introduction}

In tropical geometry, algebraic curves are replaced by certain balanced piece-wise linear graphs called tropical curves. Tropical geometry has gained lots of attention recently. One of the interesting results is that we can determine numbers of algebraic curves on toric surfaces satisfying incidence conditions by counting the corresponding tropical curves instead (Mikhalkin's Correspondence Theorem, see \cite{Mi03}).
This is true in particular for the toric surfaces $\F_0=\P^1\times \P^1$ and $\F_2$.

Gromov-Witten invariants can be thought of as ``virtual'' solutions to enumerative problems.
They are deformation invariants, thus they coincide for the two surfaces $\F_0$ and $\F_2$. For $\F_0$, Gromov-Witten invariants are enumerative, i.e.\ they count curves on $\F_0$ satisfying incidence conditions. For $\F_2$, they are not, but it is known how they are related to enumerative numbers. Therefore there is an equation relating the enumerative numbers of $\F_0$ and $\F_2$.

The purpose of this paper is to give a tropical proof of this equation for the case of rational curves, using Mikhalkin's Correspondence Theorem.

Let us introduce this equation in more details. Let $C$ denote the class of a section of $\F_2$ with $C.C=2$ and $F$ the class of the fiber of ruling.
Then the Picard group of $\F_2$ is generated by $C$ and $F$. The exceptional curve is linearly equivalent to $C-2F$. The Picard group of $\F_0$ is generated by two fibers of ruling which we will denote by $C$ and $F$ as well.
We can degenerate $\F_0$ to $\F_2$ such that the class $aC+(a+b)F$ on $\F_0$ becomes the class $aC+bF$ on $\F_2$.
Then for nonnegative $a$, $b$ with $a+b\geq 1$ we have
 
\begin{equation}
\label{eq1}
\numbc{$g$}{$0$}{$a$}{$a+b$} = \sum_{k=0}^{a-1}{\binom{b+2k}{k}\numbc{$g$}{$2$}{$a-k$}{$b+2k$}}.
\end{equation}
where $\numbc{$g$}{$0$}{$a$}{$a+b$}$ and $\numbc{$g$}{$2$}{$a$}{$b$}$ denote the numbers of nodal irreducible curves of genus $g$ of class $aC+(a+b)F$ in $\F_0$ (resp.\ of class $aC+bF$ in $\F_2$) through $4a+2b+g-1$ points in general position. (See \cite{AB01}, theorem 3.1.1, for rational curves and \cite{Vak00b}, section 8.3, for arbitrary genus.)

Let us now introduce the analogous tropical numbers.
The polygon corresponding to the divisor class $aC+(a+b)F$ on $\F_0$ is a rectangle with vertices $(0,0)$, $(a,0)$, $(a,a+b)$ and $(0,a+b)$, the polygon corresponding to the divisor class  $aC+bF$ on $\F_{2}$ is a quadrangle with vertices $(0,0)$, $(a,0)$, $(a,b)$ and $(0,2a+b)$.

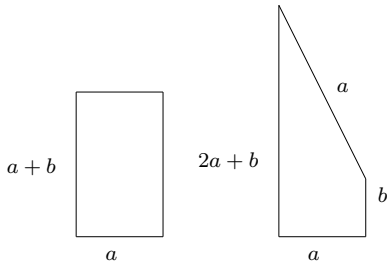
\begin{center}
\input{polytope.pstex_t}
\end{center}

We consider plane tropical curves dual to these polygons. Thus, we consider plane tropical curves of degree $\Delta_{\F_{0}}(a,a+b)$ and $\Delta_{\F_{2}}(a,b)$, where $\Delta_{\F_{0}}(a,a+b)$ denotes the multiset of the vectors $(-1,0)$ and $(1,0)$ each $a+b$ times and $(0,-1)$ and $(0,1)$ each $a$ times and $\Delta_{\F_{2}}(a,b)$ denotes the multiset of the vectors $(-1,0)$ $2a+b$ times, $(0,-1)$ $a$ times, $(1,0)$ $b$ times and $(2,1)$ $a$ times.
We denote by $\numb{$g$}{$0$}{$a$}{$a+b$}$ (resp.\ $\numb{$g$}{$2$}{$a$}{$b$}$) the number of irreducible plane \emph{tropical} curves of degree $\Delta_{\F_{0}}(a,a+b)$ (resp.\ $\Delta_{\F_{2}}(a,b)$) and genus $g$ through $4a+2b+g-1$ points in general position (see \cite{Mi03}).

Our central result is the following theorem:
\begin{theorem}\label{thm-main}
The following equation holds for nonnegative integers $a$, $b$ with $a+b\geq 1$:
\begin{equation}
\label{eq1trop}
\numb{$0$}{$0$}{$a$}{$a+b$} = \sum_{k=0}^{a-1}{\binom{b+2k}{k}\numb{$0$}{$2$}{$a-k$}{$b+2k$}}.
\end{equation}
\end{theorem}
Of course this theorem (and even the more general case) is an immediate consequence of equation \ref{eq1} and Mikhalkin's Correspondence theorem which states that $\numbc{$g$}{$0$}{$a$}{$a+b$}=\numb{$g$}{$0$}{$a$}{$a+b$}$ and $\numbc{$g$}{$2$}{$a$}{$b$}=\numb{$g$}{$2$}{$a$}{$b$}$.
However, we want to give a proof within tropical geometry.

We use induction on the degree and generalizations of Kontsevich's formula for enumerative numbers on $\F_0$ and $\F_2$.
While Kontsevich's formula for $\F_0$ (see theorem \ref{kontsevich0}) was known and can be proved without tropical geometry (\cite{FP97}, section 9), our formula for $\F_2$ (see theorem \ref{kontsevich2}) is new and was derived using tropical geometry. 
To derive a Kontsevich-type formula tropically, we compute numbers of curves satisfying point and line conditions and mapping to a special point in tropical $\mathcal{M}_{0,4}$ under the forgetful map. To prove Kontsevich's formula for $\P_2$, one can show that all such curves have a contracted bounded edge and can thus be interpreted as reducible tropical curves (\cite{GM053}). For $\F_2$, this statement is no longer true. Instead, we get a correction-term counting curves that do not have a contracted bounded edge. We show that these curves can also be interpreted as reducible curves in a different way.

Unfortunately, it seems that our method cannot be generalized to higher genus.

The paper is organized as follows. In section \ref{sec-1}, we prove our tropical Kontsevich formulas for $\F_0$ and $\F_2$. In section \ref{sec-2}, we use those formulas to prove theorem \ref{thm-main} using induction.

We would like to thank Ionut Ciocan-Fontanine and Andreas Gathmann for helpful discussions.

\section{Tropical Kontsevich formulas for $\F_0$ and $\F_2$}
\label{sec-1}
To derive tropical Kontsevich formulas for $\F_0$ and $\F_2$, we generalize the ideas of \cite{GM053}. Let us start by recalling some notations we will use. 
\begin{notation}\label{not}
Let $\Delta=\Delta_{\F_{2}}(a,b)$ and let $\M_{0,n}^{\text{lab}}(\R^2,\Delta)$ denote the space of rational parametrized tropical curves in $\R^2$ of degree $\Delta$, with $\#\Delta+n$ ends all of which are labelled, and $n$ of which are contracted ends, also called marked points (see definition 4.1 of \cite{GKM07}).
The following picture shows an example, an element $(\Gamma,h,x_1,\ldots,x_5)$ of $\M_{0,5}^{\text{lab}}(\R^2,\Delta)$ with $\Delta=\Delta_{\F_{2}}(1,1)$. On the left, we have a rational abstract tropical curve, i.e.\ a tree with unbounded edges (ends) whose vertices must have valence at least $3$. The bounded edges $e$ are equipped with lengths $l(e)\in \R_{>0}$ (which is not indicated in the picture). The map $h$ is affine-linear with rational slope on each edge. It stretches each edge of length $l(e)$ to a line segment of length $l(e)\cdot \omega(e)\cdot \Vert u(e)\Vert$, where $\omega(e)\in \N$ is called the weight of the edge, $u(e)$ is the primitive integer vector in the direction of the edge (only defined up to sign), and $\Vert . \Vert$ denotes the usual euclidean length. The vector $v(e):=\omega(e)\cdot u(e)$ is called the direction vector of the edge $e$. The map $h$ has to satisfy the balancing condition at each vertex, i.e.\ the sum of all direction vectors of edges adjacent to the vertex has to be zero. 
Contracted ends --- drawn with dotted lines and marked by $x_1,\ldots,x_5$ in the picture --- have direction vector $0$. Thus they are mapped to a point in $\R^2$, which justifies to call them marked points. The labels of the other ends are neglected in the picture. The other ends have to be mapped to the $6$ directions prescribed by $\Delta$.
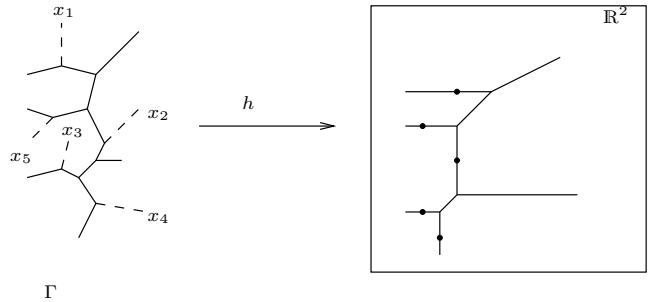
\begin{center}
\input{curvebsp.pstex_t}
\end{center}

The forgetful map $\ft$ denotes the map \begin{displaymath}\ft:\M_{0,n}^{\text{lab}}(\R^2,\Delta)\rightarrow \M_{0,4}\end{displaymath}  which forgets all ends but the first $4$ marked points (see definition 4.1 of \cite{GM053}).
This means we take the smallest subgraph of $\Gamma$ that contains $x_1,\ldots,x_4$ and straighten $2$-valent vertices. In the example, we thus have to add the two lengths $l_1$ and $l_2$ of the two bounded edges which become one bounded edge after straightening.
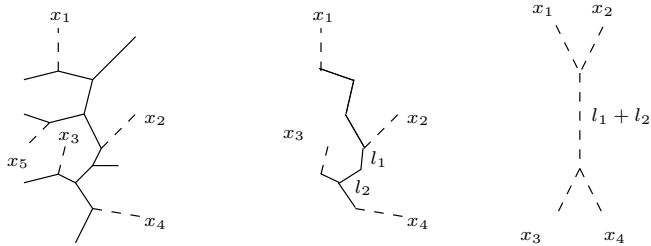
\begin{center}
\input{curvebsp2.pstex_t}
\end{center}
The evaluation map at the marked point labelled $i$ (see definition 3.3 of \cite{GM053}) is denoted by $\ev_i$:
 \begin{displaymath}\ev_i:\M_{0,n}^{\text{lab}}(\R^2,\Delta)\rightarrow \R^2: (\Gamma,h,x_i)\mapsto h(x_i).\end{displaymath} 
Pick two tropical rational functions $\varphi_{A}$ and $\varphi_{B}$ on $\M_{0,4}$ (in the sense of \cite{AR07}, definition 3.1, i.e.\ continuous piece-wise linear functions with rational slopes) whose divisors correspond to abstract tropical curves $\lambda_A$ (resp.\ $\lambda_B$) where the ends $x_1$ and $x_2$ come together at a vertex (resp.\ where $x_1$ and $x_3$ come together) and where the length parameter of the bounded edge is very large.
\end{notation}

\begin{remark}
Note that we use the space of parametrized tropical curves with labelled ends here. The reason is that one can show that this space is a tropical fan (proposition 4.7 of \cite{GKM07}) and that we can thus use the tropical intersection theory from \cite{AR07}. Since we want to count tropical curves without the labels of the non-contracted ends, we have to divide by a factor of $|G|$, where $G$ is the group of possible permutations of the labels.
In the Kontsevich formula we want to prove (theorem \ref{kontsevich2}), we sum over all possibilities to split the degree $\Delta=\Delta_{\F_{2}}(a,b)$ into two smaller degrees.
To be precise, we would have to sum over all possibilities to pick a \emph{labelled} subset of non-contracted ends forming the smaller degrees. This factor together with the factors for labelling the ends in the small degrees exactly cancel with the total factor of $|G|$. In the following, we will therefore neglect the fact that non-contracted ends are labelled. 
\end{remark}

 The difference of $\varphi_{A}$ and $\varphi_{B}$ is globally given by a bounded tropical rational function on $\M_{0,4}$. Therefore, the tropical Cartier divisors $[\varphi_{A}]$ and $[\varphi_{B}]$ are rationally equivalent and by lemma 8.5 of \cite{AR07} their pull-backs  $[\ft^{\ast}\varphi_{A}]\cdot \M_{0,n}^{\text{lab}}(\R^2,\Delta)$ and $[\ft^{\ast}\varphi_{B}]\cdot\M_{0,n}^{\text{lab}}(\R^2,\Delta)$ are rationally equivalent as well.

Set $n=\# \Delta$ and choose tropical rational functions $\varphi_{1}$, $\varphi_{2}$, ${\varphi_{3}}_{1}, {\varphi_{3}}_{2}$ $\ldots$, ${\varphi_{n}}_{1}, {\varphi_{n}}_{2}$ on $\R^{2}$ that correspond to tropical curves $L_{1}$ and $L_{2}$  of degree $\Delta_{\F_{2}}(1,0)$ and to points $p_{3}, \ldots, p_{n} \in \R^{2}$ in general position. 
 We can set $\varphi_{1}=\max\{x-{p_{1}}_{1},2(y-{p_{1}}_{2}),0\} $ and $\varphi_{2}=\max\{x-{p_{2}}_{1},2(y-{p_{2}}_{2}),0\}$ to get $L_1$ and $L_2$ in this case.

Because of the above, we have 
\begin{equation}\label{eqdeg}
\begin{split}
&\deg([\ev_{1}^{\ast}\varphi_{1}\cdot\ev_{2}^{\ast}\varphi_{2}\cdot\prod_{i=3}^n(\ev_{i}^{\ast}{\varphi_{i}}_{1}\cdot\ev_{i}^{\ast}{\varphi_{i}}_{2})\cdot\ft^{\ast}\varphi_{A}\cdot\M_{0,n}^{\text{lab}}(\R^2,\Delta)]) \\
=&\deg([\ev_{1}^{\ast}\varphi_{1}\cdot\ev_{2}^{\ast}\varphi_{2}\cdot\prod_{i=3}^n(\ev_{i}^{\ast}{\varphi_{i}}_{1}\cdot\ev_{i}^{\ast}{\varphi_{i}}_{2})\cdot\ft^{\ast}\varphi_{B}\cdot\M_{0,n}^{\text{lab}}(\R^2,\Delta))].
\end{split}
\end{equation} 
\begin{remark}\label{rem-generalconditions}
Both above expressions are $0$-dimensional tropical intersection products as defined in \cite{AR07}, even if the set-theoretical intersection is higher-dimensional. If we pick the conditions to be general however, the set-theoretical intersection equals the support of the intersection product. That means that the intersection products above equal the sums of tropical curves in $\M_{0,n}^{\text{lab}}(\R^2,\Delta)$ that satisfy the conditions, i.e.\ that pass through $L_1$, $L_2$, $p_3,\ldots,p_n$ and map to $\lambda_A$ resp.\ $\lambda_{B} \in \M_{0,4}$ under $\ft$, counted with multiplicity. This can be shown analogously to \cite{MR08}, lemma 3.1.
\end{remark}
The following lemma will enable us to compute the multiplicity with which we have to count each curve satisfying the conditions in the intersection product:

\begin{lemma}
\label{lem-determinante}
Let $X$ be an abstract tropical cycle (in the sense of \cite{AR07}, definition 5.12, i.e.\ roughly a weighted polyhedral complex satisfying the balancing condition) of dimension $k$ and $\varphi_{1}$, \ldots, $\varphi_{k}$ tropical rational functions on $X$. Moreover, let $P \in X$ be a point in the interior of a cone $\sigma$ of maximal dimension in $X$. Assume that $\varphi_{i}$ is of the form $\max\{ \psi_{i} , \chi_{i} \}$ locally around $P$, where $\psi_{i}, \chi_{i}: X \rightarrow \R$ denote $\Z$-affine functions with $\psi_{i}(P)=\chi_{i}(P)$. 
Let $(\psi_{i}-\chi_{i})_L$ denote the linear part of the affine function $(\psi_{i}-\chi_{i})$ and let $A$ be the $(k \times k)$-matrix with entries $((\psi_{i}-\chi_{i})_L(u_{j}))_{i,j}$ for a basis $u_{1},  \ldots, u_{k}$ of the lattice underlying $X$ at $\sigma$. Then the coefficient of $P$ in the intersection product $\varphi_{1}\cdot\ldots\cdot\varphi_{k}\cdot X$ is equal to $\omega(\sigma)\cdot |\det(A)|$.
\end{lemma}

\begin{proof}
The computation of the coefficient of $P$ in the intersection product is local around $P$. Thus, we may assume that $X$ is a tropical fan in some vector space $V$ and extend $\sigma$ to the affine vector space $V_{\sigma}$ spanned by $\sigma$. Furthermore, we may consider the tropical rational functions $\varphi_{i}=\max\{\psi_{i}, \chi_{i}\}$ on the whole space $V_{\sigma}$. Moreover, we may replace the tropical rational functions $\max\{\psi_{i},\chi_{i}\}$ by $\max\{\psi_{i},\chi_{i}\}-\chi_{i}=\max\{\psi_{i}-\chi_{i},0\}$ as changing a tropical rational function by a linear function does not affect the intersection product. We define the morphism $g = (\psi_{1}-\chi_{1}, \ldots, \psi_{k}-\chi_{k}): X \rightarrow \R^{k}$. Then we have $\varphi_{i}=g^{\ast}\mu_{i}$ for $\mu_{i}:\R^{k} \rightarrow \R; (a_{1}, \ldots, a_{k}) \mapsto \max\{a_{i},0\}$ and for all $1 \leq i \leq k$. By the projection formula (\cite{AR07}, proposition 4.8) the multiplicity of $P \in X$ in the intersection product $\varphi_{1} \cdot \ldots \cdot \varphi_{k} \cdot \sigma = g^{\ast}\mu_{1} \cdot \ldots \cdot g^{\ast}\mu_{k} \cdot \sigma$ is equal to the multiplicity of $0$ in $\R^{k}$ in the intersection product $g_{\ast}(g^{\ast}\mu_{1} \cdot \ldots \cdot g^{\ast}\mu_{k} \cdot \sigma) = \mu_{1} \cdot \ldots \mu_{k} \cdot g_{\ast}\sigma$. For dimensional reasons the cycle $g_{\ast}\sigma$ is the whole target space $\R^{k}$ with some weight. But this weight is $\omega(\sigma)\cdot|\det(A)|$. Note $\mu_{1} \cdot \ldots \cdot \mu_{k} \cdot \R^{k}$ is the origin with weight $1$. This finishes the proof.
\end{proof}

\begin{remark}\label{rem-choosecoor}
If $\sigma$ is a cone in $\mathcal{M}^{\text{lab}}_{0,n}(\R^2,\Delta)$, it corresponds to a combinatorial type, i.e.\ a homeomorphism class of a graph plus direction vectors for all edges (see \cite{GM053}, 2.9). We can deform a parametrized tropical curve $(\Gamma,h,x_i)$ within $\sigma$ by changing the length of the bounded edges or translating the image $h(\Gamma)$.
  Thus a basis for the lattice underlying $\mathcal{M}^{\text{lab}}_{0,n}(\R^2,\Delta)$ at $\sigma$ is given by the position of a root vertex $h(V)$ and the length of all bounded edges. By remark 3.2 of \cite{GM053}, the absolute value of the determinant of the matrix $A$ from lemma \ref{lem-determinante} above is independent from the choice of a root vertex and an order of the bounded edges.
\end{remark}

\begin{example}\label{ex-det}
Assume $\sigma$ is the cone in $\mathcal{M}^{\text{lab}}_{0,4}(\R^2,\Delta_{\F_2}(1,0))$ corresponding to the combinatorial type pictured below.
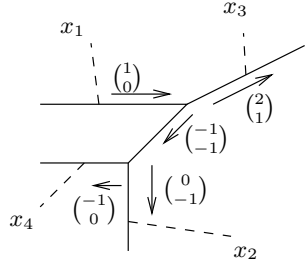
\begin{center}
\input{detbspcomb.pstex_t}
\end{center}
Following remark \ref{rem-choosecoor}, we choose the position of $h(x_1)$ and the lengths $l_1,\ldots,l_5$ of the bounded edges as coordinates for $\sigma$.
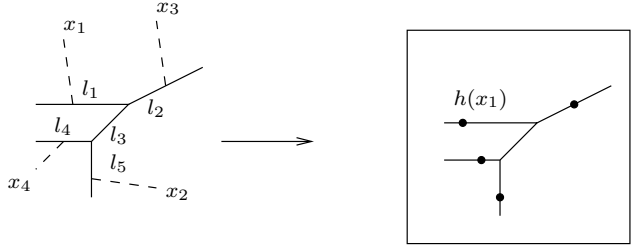
\begin{center}
\input{detbspcurve.pstex_t}
\end{center}
Consider the following curve $C$ inside $\sigma$:

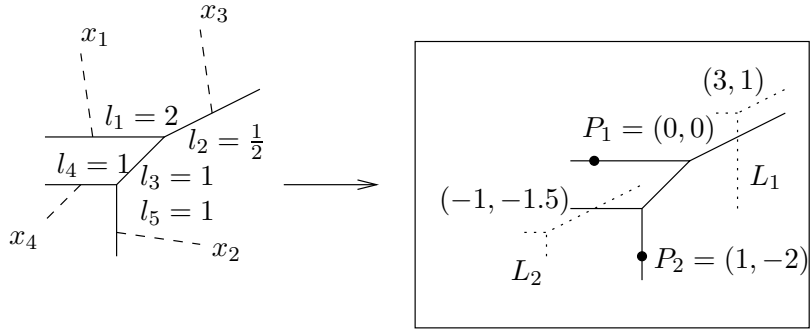
\begin{center}
\input{detbsp.pstex_t}
\end{center}

The curve $C$ (where $h(x_1)=(0,0)$, $l_1=2$, $l_2=\frac{1}{2}$ and $l_3=l_4=l_5=1$) goes through the points $P_1$ (which is cut out by $\max\{x,0\}$ and $\max\{y,0\}$) and $P_2$ (cut out by $\max\{x,1\}$ and $\max\{y,-2\}$) and through $L_1$ (cut out by $\max\{x-3,2y-2,0\}$) and $L_2$ (cut out by $\max\{x+1,2y+3,0\}$) and maps to the abstract tropical curve with $x_1$ and $x_3$ at one vertex and length parameter $1$ under $\ft$. Denote by $\lambda$ a tropical rational function on $\mathcal{M}_{0,4}$ that cuts out this curve.
Because $C$ satisfies the conditions, it contributes to the intersection product
\begin{align*}&\ev_1^\ast(\max\{x,0\})\cdot\ev_1^\ast(\max\{y,0\}) \ev_2 ^\ast(\max\{x,1\})\cdot \ev_2^\ast(\max\{y,-2\})\\ \cdot & \ev_3^\ast(\max\{x-3,2y-2,0\})\cdot \ev_4^\ast(\max\{x+1,2y+3,0\})\\\cdot& \ft^\ast(\lambda)\cdot\mathcal{M}_{0,4}^{\text{lab}}(\R^2,\Delta_{\F_2}(1,0)). \end{align*}
Let us compute the multiplicity with which it contributes using lemma \ref{lem-determinante}.
Locally at $h(x_3)$, the function $\max\{x-3,2y-2,0\}$ equals $\max\{x-3,0\}$ and locally at $h(x_4)$, the function $\max\{x+1,2y+3,0\}$ equals $\max\{x+1,2y+3\}$. Hence locally we have \begin{align*}\ev_3^\ast(\max\{x-3,2y-2,0\})&=\max\{h(x_3)_x-3,0\} \mbox{ and }\\ \ev_4^\ast(\max\{x+1,2y+3,0\})&=\max\{h(x_4)_x+1,2h(x_4)_y+3\}.\end{align*} We can rewrite the pullbacks along $\ev_1$ and $\ev_2$ analogously. Locally, $\ft$ equals the map that sends a curve in $\sigma$ with coordinates $(h(x_1),l_1,\ldots,l_5)$ to $l_3$. We also have to write the linear part of the evaluation pullbacks in the basis of $\sigma$, i.e. in $h(x_1),l_1,\ldots,l_5$.
For this, note first that
\begin{align*}
h(x_2)&=h(x_1)+l_1\cdot \binom{1}{0}+l_3\cdot\binom{-1}{-1}+l_5\cdot\binom{0}{-1},\\
h(x_3)&=h(x_1)+l_1\cdot\binom{1}{0}+l_2\cdot\binom{2}{1}\mbox{ and }\\
h(x_4)&=h(x_1)+l_1\cdot \binom{1}{0}+l_3\cdot\binom{-1}{-1}+l_4\binom{-1}{0}.
\end{align*}
Thus, the linear part of $h(x_3)_x-3$ equals \begin{displaymath}h(x_1)_x+l_1+2l_2\end{displaymath} and the linear part of $h(x_4)_x+1-2h(x_4)_y-3$ equals \begin{displaymath}h(x_1)_x+l_1-l_3-l_4-2(h(x_1)_y-l_3)=h(x_1)_x-2h(x_1)_y+l_1+l_3-l_4.                    \end{displaymath}
So, if we plug for example the vector which has a $1$ at $l_2$ and $0$ everywhere else into the linear part of $h(x_3)_x-3$, we get $2$.
If we plug the vector which has a $1$ at $l_4$ and $0$ everywhere else into the linear part of $h(x_4)_x+1-2h(x_4)_y-3$, we get $-1$. Continuing like this, we can see that the matrix $A$ equals:
 \[ \begin {array}{l|ccccccc}
&h(x_1)_x&h(x_1)_y& l_1&l_2&l_3&l_4&l_5\\ \hline
(P_1)_x &1&0&0&0&0&0&0\\
(P_1)_y &0&1&0&0&0&0&0\\
(P_2)_x &1&0&1&0&-1&0&0\\
(P_2)_y &0&1&0&0&-1&0&-1\\
L_1&1&0&1&2&0&0&0\\
L_2&1&-2&1&0&1&-1&0\\
\ft &0&0&0&0&1&0&0
         \end {array} \]

Since $ |\det(A)|=2$, the curve $C$ contributes with multiplicity $2$ to the intersection product above.

\end{example}

\begin{notation}\label{def-mult}
Let $C$ be a curve contributing to a $0$-dimensional intersection product as in example \ref{ex-det} consisting of evaluation pullbacks and the pullback of a curve in $\mathcal{M}_{0,4}$ under $\ft$ (resp.\ only evaluation pullbacks).
Then we denote by $\mult_{\ev\times \ft}(C)$ (resp.\ $\mult_{\ev}(C)$) the multiplicity with which $C$ contributes to the intersection product, which equals the absolute value of the determinant of the linear part of the combined evaluation and forgetful maps, as we have seen in \ref{lem-determinante}.
\end{notation}

\begin{remark}\label{rem-detmult}
Note that by \cite{GM053}, proposition 3.8, $\mult_{\ev}(C)$ equals the usual multiplicity of a tropical curve as defined in \cite{Mi03}, 4.15, i.e.\ the multiplicity with which it contributes to the count of $\numb{$g$}{$2$}{$a$}{$b$}$.
\end{remark}

In the following, we want to describe both sides of equation \ref{eqdeg} in detail. We want to study the set of curves that satisfy the conditions, and their multiplicity. We will see that we can interpret the curves as reducible curves, and count the contributions from each component separately. This will lead to the formula of theorem \ref{kontsevich2} we want to prove.

\begin{remark}\label{rem-string}
Using the notations from \ref{not}, let $C$ be a tropical curve in $\M_{0,n}^{\text{lab}}(\R^2,\Delta)$ passing through $L_1$, $L_2$, $p_3,\ldots,p_n$ and mapping to $\lambda_A$ under $\ft$ (hence a curve $C$ that contributes to the left hand side of equation \ref{eqdeg} with multiplicity $\mult_{\ev\times \ft}(C)$).
We would like to generalize proposition 5.1 of \cite{GM053}, which states that $C$ has a contracted bounded edge. However, this is not true in the case of $\F_2$. We can have curves like the one shown in the following picture (where the lengths $l$ and $l'$ are very large) which do allow a very large $\mathcal{M}_{0,4}$-coordinate.
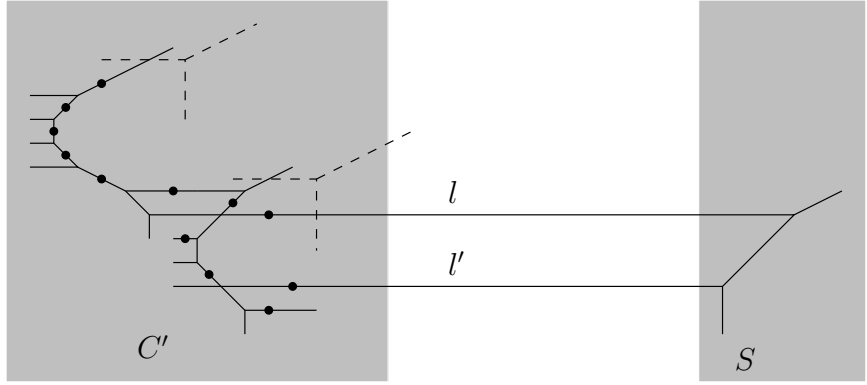
\begin{center}
\input{largem04.pstex_t}
\end{center}
Even though those curves fail to have a contracted bounded edge, we can still interpret them as reducible curves by cutting off the part which is far away to the right (in the picture denoted by $S$). The remaining part (in the picture denoted by $C'$) is a reducible curve of degree $\Delta_{\F_2}(a-1,b+2)$. The existence of such curves with a very large $\mathcal{M}_{0,4}$-coordinate leads to the second part of the sum in the recursion formula of theorem \ref{kontsevich2}. The part $S$ which is far away to the right is called a \emph{string} following \cite{GM053}, definition 3.5, i.e.\ a subgraph of the abstract tropical curve $\Gamma$ homeomorphic to $\R$ (a ``path in $\Gamma$ with two unbounded ends'') that does not intersect the closures $\overline{x_i}$ of the marked points. The two ends are of direction $(0,-1)$ and $(2,1)$ as in the picture.
\end{remark}

\begin{lemma}\label{lem-contractededge}
Using the notations from \ref{not}, let $C \in \M_{0,n}^{\text{lab}}(\R^2,\Delta) $ be a tropical curve that passes through $L_{1}$, $L_{2}$, $p_{3}$, \ldots, $p_{n}$, maps to $\lambda_A$ under $\ft$ and has a non-zero multiplicity $\mult_{\ev\times \ft}(C)$. Then \textbf{either}
\begin{enumerate}
\item $C$ has a contracted bounded edge \textbf{or}
\item $C$ contains a string $S$ (see remark \ref{rem-string}) that can be moved to the right.
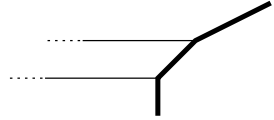
\label{case-string}
\begin{center}
\input{stueck.pstex_t}
\end{center}
\end{enumerate}
\end{lemma}

\begin{proof}
The beginning of the proof is similar to proposition 5.1 of \cite{GM053}.

We will show that the set of all points $\ft(C)$ is bounded in $\M_{0,4}$ where $C$ runs over all curves $C \in \M_{0,n}^{\text{lab}}(\R^2,\Delta)$ with non-zero multiplicities $\mult_{\ev\times \ft}(C)$ that satisfy the conditions but have no contracted bounded edge and no string $S$ with ends $(0,-1)$ and $(2,1)$ moving to the right as in the picture. By proposition 2.11 of \cite{NS06} there are only finitely many combinatorial types in $\M_{0,n}^{\text{lab}}(\R^2,\Delta) $. Thus, we may restrict ourselves to tropical curves $C $ of a fixed combinatorial type $\alpha$. Furthermore, we may assume the curves corresponding to $\alpha$ are $3$-valent.

Let $C $ be such a curve and let $C'$ be the curve obtained from $C$ by forgetting the first and the second marked point. Then $C'$ has a string $\Gamma'$ which follows analogously to remark 3.7 of \cite{GM053} by counting connected components of the abstract graph without the closures of the marked points --- there are less such components than ends. As explained below, we can use the string to deform $C'$ in a $1$-parameter family within its combinatorial type without changing the images of the marked points.
Each such deformation of $C'$ yields a curve satisfying the point conditions and meeting the two curves $L_1$ and $L_2$. We show that the $1$-dimensional deformation of $C'$ is either bounded itself or does not affect the image under $\ft$. From this, the statement follows.

Using the string, we can deform $C'$: Starting with one of the ends of $\Gamma'$ and in a nonzero direction modulo its direction vector, we can translate the images $h(e)$ of edges $e$ of the string $\Gamma'$, prolonging or shortening some of the other edges, but keeping the lines in $\R^2$ fixed to which these other edges are mapped.
Since we want to stay within the combinatorial type however, we cannot shorten any length to $0$, thus shrinking an edge.
This means that if there is a bounded edge $e$ which is mapped e.g.\ to the left of $h(\Gamma')$ (after picking an orientation of $\Gamma'$), then the translation of $h(\Gamma')$ to the left is bounded by the length of this edge (as in (i) of the picture below).
In particular, if there are bounded edges on both sides of $h(\Gamma')$, the deformation of $C'$ with the combinatorial type and the conditions fixed are bounded. This means that the lengths of all inner edges are bounded except possibly the edges adjacent to $x_{1}$ and $x_{2}$. This is sufficient to ensure that the image of these  curves under $\ft$ is bounded in $\M_{0,4}$, too. Hence we can assume now that we have bounded edges only on one side of $h(\Gamma')$.
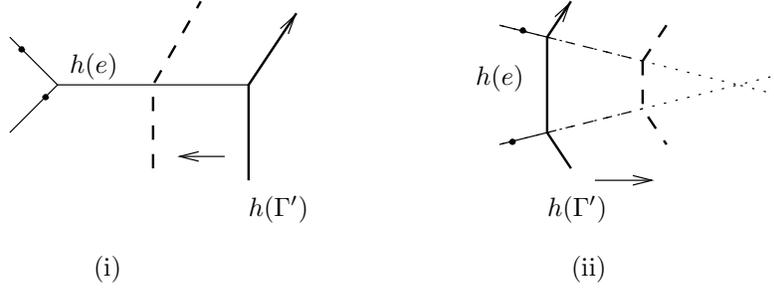
\begin{center}
\input{stringmov1.pstex_t} 
\end{center}

Assume we could deform $C'$ in a more than $1$-dimensional family while keeping the images of the marked points fixed. Then $C$ moves in an at least $1$-dimensional family with the image point
  under all evaluations and the forgetful map fixed. Then $\mult_{\ev\times \ft}(C)= 0$, which is a contradiction to our assumption. In particular, we can see that we cannot have more than one string.
Note that the edges adjacent to $\Gamma'$ must be bounded since otherwise we would have two strings.
So now we can assume that there are only bounded edges adjacent to the string, and that they are all on one side of $h(\Gamma')$ (say after picking an orientation of $\Gamma'$ on the left side). Denote the direction vectors of the edges of $\Gamma'$ by $v_{1}, \ldots, v_{k}$ and the direction vectors of the adjacent bounded edges by $w_{1}, \ldots, w_{k-1}$. We have seen already that we cannot translate $h(\Gamma')$ without limit to the left. If one of the directions $w_{i+1}$ is obtained from $w_{i}$ by a right turn, then the lines to which the edges corresponding to $w_{i}$ and $w_{i+1}$ are mapped meet to the right of $\Gamma'$ as shown in (ii). This restricts the movement of $\Gamma'$ to the right with the combinatorial type and the conditions fixed, too, since the edge corresponding to $v_{i+1}$ then receives length $0$.
Hence, as above, the image of these plane tropical curves under $\ft$ is bounded in $\M_{0,4}$ as well. Thus, we may assume that for all $i$, $1 \leq i \leq k-2$ the direction $w_{i+1}$ is either the same as $w_{i}$ or obtained from $w_{i}$ by a left turn as shown in (iii) of the picture below. The balancing condition then ensures that for all $i$ both the directions $v_{i+1}$ and $-w_{i+1}$ lie in the angle between $v_{i}$ and $-w_{i}$. Therefore, all directions $v_{i}$ and $-w_{i}$ lie in the angle between $v_{1}$ and $-w_{1}$. In particular, the string $\Gamma'$ cannot have any self-intersections in $\R^{2}$. We can therefore pass to the local dual picture where the edges dual to $w_{i}$ correspond to a concave side of a polygon whose other two edges are dual to $v_{1}$ and $v_{k}$ as shown in (iv).
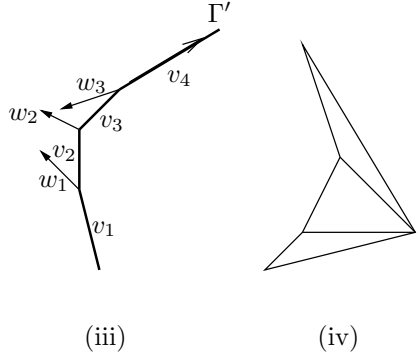
\begin{center}
\input{contracted2.pstex_t}
\end{center}
But note that both $v_{1}$ and $v_{k}$ are outer directions of a plane tropical curve of degree $\Delta$. 
 Thus, $v_{1}$ and $v_{k}$ must be $\binom{-1}{0}$, $\binom{0}{-1}$, $\binom{1}{0}$ or $\binom{2}{1}$. Consequently, their dual edges have direction vectors $\pm\binom{0}{-1}$, $\pm\binom{1}{0}$, $\pm\binom{0}{1}$ or $\pm\binom{-1}{2}$. We have to distinguish two cases
\begin{itemize}
\item[(a)] $v_{1}$ and $v_{k}$ are $\binom{0}{-1}$ and $\binom{2}{1}$, i.e.\ their dual edges have direction vectors $\pm(-1,0)$ and $\pm(-1,2)$
\item[(b)] $v_{1}$ and $v_{k}$ are not $\binom{0}{-1}$ and $\binom{2}{1}$.
\end{itemize}

In case (b) the triangles spanned by two of those vectors do not
admit any further integer points. 
Therefore we have $k=2$ and the string consist just of the
two unbounded edges corresponding to $v_{1}$ and $v_{2}$ that are
connected to the rest of the plane tropical curve by exactly one
internal edge corresponding to $w_{1}$. 
We now do a case-by-case analysis going through all $5$ choices of $v_1$ and $v_k$ which are not $(0,-1)$ and $(2,1)$. In each case, the direction $w_1$ is fixed. As the length of the edge corresponding to $w_1$ increases, the image of the union of the two ends corresponding to $v_1$ and $v_k$ becomes disjoint from at least one of the chosen curves $L_1$ and $L_2$. Hence the two marked points $x_1$ and $x_2$ cannot be adjacent to those two ends, or we again leave the combinatorial type with the deformation. But if $x_1$ and $x_2$ are not adjacent to those two ends, the lengths of the edge corresponding to $w_1$ does not affect the image under $\ft$ and thus again, the image of the $1$-dimensional family under $\ft$ is bounded. 
Below, we depict $4$ of the $5$ cases we have to consider. The fifth case is not possible in fact, since it would lead to a contracted bounded edge and thus to a more than $1$-dimensional deformation.

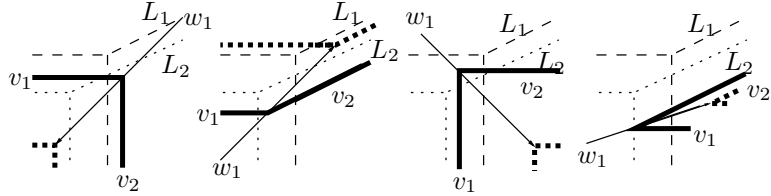
\begin{center}
\input{contractedf2.pstex_t}
\end{center}

In case (a) the triangle spanned by the two vectors $\binom{-1}{0}$ and $\binom{-1}{2}$ admits exactly one further integer point.

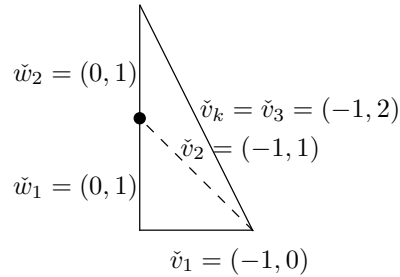
\begin{center}
\input{nocontractedA.pstex_t}
\end{center}

In the picture, we denote the duals of the vectors $v_i$ and $w_i$ by $\check{v}_i$ and $\check{w}_i$.
Thus, in case (a) we may have $k=3$ and the string $\Gamma'$ may consist of the two unbounded edges corresponding to $v_{1}$ and $v_{3}$ and the bounded edge corresponding to $v_{2}$ that is connected to the rest of the plane tropical curve by the two edges corresponding to $w_{1}$ and $w_{2}$. In this case, the movement of the string is indeed not bounded to the right.
Then we are in case (\ref{case-string}) of lemma \ref{lem-contractededge}. This finishes the proof of the lemma.
\end{proof}

\begin{lemma}\label{lem-multstring}
Let $C$ be a curve of type (\ref{case-string}) of lemma \ref{lem-contractededge} then 
\[ \mult_{\ev\times \ft}(C) = \mult_{\ev_1}(C_1) \cdot \mult_{\ev_2}(C_2) \cdot 2 \cdot (C_{1} \cdot L_{1})_{x_{1}} \cdot (C_{1} \cdot L_{2})_{x_{2}} \]
where $\mult_{\ev_{i}}(C_i)$ denotes the multiplicity of the evaluation map at the $\#\Delta_{i}-1$ points of $x_{3}, \ldots, x_{n}$ that lie on $C_{i}$ for $i \in \{1,2\}$ and $(C' \cdot C'')_{p}$ denotes the intersection multiplicity of the plane tropical curves $C'$ and $C''$ at the point $p\in C'\cap C''$.
Here, $C_1$ and $C_2$ denote the two irreducible components of the part $C'$ of $C$ that we get when cutting off the string $S$ as in remark \ref{rem-string}.
\end{lemma}

\begin{proof}
Since $\mult_{\ev\times \ft}(C)$ equals the absolute value of the determinant of the map $\ev\times \ft$ in local coordinates, we set up the matrix $A$ for $\ev\times \ft$  as in lemma \ref{lem-determinante} and compute its determinant.
The local coordinates are the position of a root vertex and the length of all bounded edges, respectively the coordinates of the images of the marked points and the length coordinate of the bounded edge of the image under $\ft$ in $\M_{0,4}$. Because of remark \ref{rem-choosecoor}, the absolute value of the determinant does not depend on the special choice of such coordinates.

There are exactly two bounded edges that connect the string $S$ with the rest of the curve. We denote these bounded edges by $E'$ and $E''$ and the unique bounded edge that is contained in the string by $E$. Their lengths are denoted by $l'$, $l''$ and $l$, respectively. 

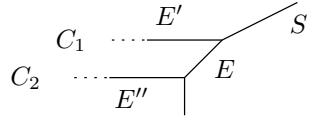
\begin{center}
\input{stringS.pstex_t}
\end{center}

As the length of the $\M_{0,4}$-coordinate is very large and there is no contracted bounded edge, the lengths $l'$, $l$ and $l''$ must count towards the length of the $\M_{0,4}$-coordinate. 
That is, $x_1$ and $x_2$ have to be on one side of those three edges and $x_3$ and $x_4$ on the other. Let us call the part with $x_1$ and $x_2$ $C_1$, and assume without restriction that $E'$ belongs to $C_1$.
Let the root vertex be a vertex in $C_1$. To evaluate at a marked point in $C_1$, we have to follow the path in $\Gamma$ from the root vertex to this marked point. We do not cross $E'$ or $E''$ with this path, thus we get $0$ in the $l'$- and $l''$-column in the rows for marked points in $C_1$. To evaluate at a marked point in $C_2$, we need to cross $E'$ in direction $(1,0)$ and $E''$ in direction $(-1,0)$. Altogether, the columns of the matrix $A$ corresponding to the lengths $l'$ and $l''$ read:

\vspace{2ex}
\begin{tabular}{l|cc}
 & $l'$ & $l''$ \\
\hline
evaluation at points behind $E'$ & $\begin{array}{c} 0 \\ 0 \end{array}$ & $\begin{array}{c} 0 \\ 0 \end{array}$ \\
\hline
evaluation at points behind $E''$ & $\begin{array}{c} 1 \\ 0 \end{array}$ & $\begin{array}{c} -1 \\ 0 \end{array}$  \\
\hline
$\M_{0,4}$-coordinate & 1 & 1 \\
\end{tabular}
\vspace{2ex}

If we add the column corresponding to the length $l'$ to the column corresponding to the length $l''$, then the new column has only one entry $2$ and all other entries $0$. Thus, we get a factor of $2$ and to compute the determinant of the matrix $A$ we may drop both the $\M_{0,4}$-row and the column corresponding to the edge $E''$.

Now, we consider the first marked point $x_{1}$. We require that
the plane tropical curve $C$ passes through $L_{1}$ at this point. Let $E_{1}$ and
$E_{2}$ be the two adjacent edges of $x_{1}$. We denote their
common direction vector by $v=(v_{1},v_{2})$ and their lengths by
$l_{1}$ and $l_{2}$, respectively. We may assume that the root
vertex is on the $E_{1}$-side of $x_{1}$. Assume that both
$E_{1}$ and $E_{2}$ are bounded. Assume $x_{1}$ is contracted to a point on an unbounded edge of $L_{1}$ with direction vector $(u_{1},u_{2})$.
Thus locally around $h(x_1)$, $L_1$ is cut out by the function $\max\{u_2x+a,u_1y+b\}$, with some affine parts $a$ and $b$. To write the evaluation at $x_1$ in our local coordinates, we have to follow the path from the root vertex to $x_1$ and produce a sum $h(V)+\sum_e l(e) \cdot v(e)$ running over the edges $e$ in the path with the directions $v(e)$ according to the path. Somewhere in this sum the summand $l_1\cdot (v_1,v_2)$ appears.
Following lemma \ref{lem-determinante}, we now have to pull back the function with the evaluation map, and then express the linear part of the difference $u_2h(x_1)_x+a-u_1h(x_1)_y-b$ in our local coordinates.
In particular, we want to understand the $l_1$-entry of this row in our matrix, i.e.\ we plug in $l_1=1$ and $0$ for all other local coordinates. Since the evaluation $h(x_1)_x$ contains a summand of $l_1\cdot v_1$ and $h(x_1)_y$ contains a summand of $l_1\cdot v_2$, the entry in our matrix is $u_2 v_1-u_1v_2$.
Notice also that the edge $E_2$ does not appear in the path from the root vertex to $x_1$, i.e.\ we have zeroes in the $l_2$-column in the row corresponding to the evaluation at $x_1$.
For any marked point that is on the $E_1$-side of $x_1$, we do not cross any of the two edges $E_1$ and $E_2$ when evaluating, i.e.\ we get zeroes in the corresponding rows. For the marked points on the $E_2$-side, we have to cross both $E_1$ and $E_2$ in direction $v$.
Altogether, the $l_1$- and $l_2$-columns of the matrix read:

\vspace{2ex}

\begin{tabular}{l|cc}
evaluation at ... & $l_{1}$ & $l_{2}$ \\
\hline
... $x_{1}$ & $u_{2}v_{1}-u_{1}v_{2}$ & $0$  \\
... a point reached via $E_{1}$ from $x_{1}$ & $0$ & $0$ \\
... a point reached via $E_{2}$ from $x_{1}$ & $v$ & $v$ \\
\end{tabular}

\vspace{2ex}

Note that there are two rows for each marked point $x_{3}$, \ldots, $x_{n}$ that is reached via $E_{1}$ or $E_{2}$ from $x_{1}$ and there is one row for the marked point $x_{2}$.
If we subtract the column corresponding to $l_{2}$ from the column corresponding to $l_{1}$, then we obtain a column with only one non-zero entry. So for the absolute value of the determinant we get a factor of $|u_{2}v_{1}-u_{1}v_{2}|$, which equals the intersection multiplicity of $C_1$ and $L_1$ at the point $h(x_1)$.
Thus we get $(C_{1} \cdot L_{1})_{x_{1}}$ as a factor and may drop both the row corresponding to $x_{1}$ and the column corresponding to $l_{1}$. Now assume that one of the edges $E_{1}$ and $E_{2}$ is unbounded. Assume that $E_{1}$ is bounded and $E_{2}$ is unbounded. Then there is a column corresponding to $l_{1}$ but no column corresponding to $l_{2}$. The column corresponding to $l_{1}$ has only one non-zero entry and the same argument as above holds. Note that it is not possible that both $E_{1}$ and $E_{2}$ are unbounded. Taking the factor $(C_{1} \cdot L_{1})_{x_{1}}$ into account, we can now forget the marked point $x_1$ and straighten the $2$-valent vertex to produce only one edge out of $E_1$ and $E_2$.

If we forget the marked point $x_{2}$, we obtain a factor of $(C_{1} \cdot L_{2})_{x_{2}}$ in the same way.

Now, we consider again the string $S$. Remember that we split up the plane tropical curve $C$ at this string into the two parts $C_{1}$ and $C_{2}$. We choose the boundary vertex of the bounded edge $E'$ at the $C_{1}$-side as root vertex and denote it by $V$. Then the matrix reads:

\vspace{2ex}
\begin{tabular}{l|ccccc}
 & & lengths & & & lengths \\
evaluation at a ... & root & in $C_{1}$ & $l'$ & $l$ & in $C_{2}$ \\
\hline
... point behind $E'$ & $I_{2}$ & $\ast$ & $0$ & $0$ & $0$ \\
\hline
... point behind $E''$ & $I_{2}$ & $0$ & $\binom{1}{0}$ & $\binom{-1}{-1}$ & $\ast$ \\
\end{tabular}

\vspace{2ex}

where $I_2$ denotes the $2$ by $2$ unit matrix.
Assume there are $n_1$ marked points besides $x_1$ and $x_2$ on $C_1$ and $n_2$ marked points besides $x_3$ and $x_4$ on $C_2$, where $n_1+n_2=n-4=\#\Delta-4=4a+2b-4$.
Assume the degree of $C_1$ is $\Delta_{\F_2}(a_1,b_1)$ and the degree of $C_2$ is $\Delta_{\F_2}(a_2,b_2)$. 
Then $a_1+a_2=a-1$ and $b_1+b_2=b+2$ as we observed in remark \ref{rem-string}.
Since $\mult_{\ev\times \ft}(C)\neq 0$ we must have $n_1=4a_1+2b_1-1$ and $n_2+2=4a_2+2b_2-1$ (because then the curves will be fixed by the points).
Thus (after forgetting $x_1$ and $x_2$) $C_1$ has $n_1+4a_1+2b_1=2n_1+1$ unbounded edges and thus $2n_1-2$ bounded edges. Hence $2n_1-2$ length coordinates belong to bounded edges in $C_1$.
$C_2$ has $n_2+2+4a_2+2b_2=2n_2+5$ unbounded edges and thus $2n_2+2$ length coordinates belong to $C_2$.
 Furthermore, there are $n_{1}$ points behind $E'$ and there are $n_{2}+2$ points behind $E''$.

If we add the $l'$-column to the $l$-column and then multiply the $l$-column by $-1$, then we obtain the following matrix whose determinant has the same absolute value as the determinant that we are looking for:

\vspace{2ex}
\begin{tabular}{l|cc|ccc}
 & & lengths & & & lengths \\
evaluation at a ... & root & in $C_{1}$ & $l'$ & $l$ & in $C_{2}$ \\
\hline
... point behind $E'$ & $I_{2}$ & $\ast$ & $0$ & $0$ & $0$ \\
\hline
... point behind $E''$ & $I_{2}$ & $0$ & $\binom{1}{0}$ & $\binom{0}{1}$ & $\ast$ \\
\end{tabular}
\vspace{2ex}

Note that this matrix has a block form. The block at the top right is a zero block. We denote the top left block of size $2n_{1}$ by $A_{1}$ and the bottom right block of size $2n_{2}+4$ by $A_{2}$. Then, the determinant that we are looking for is $|\det(A_{1})|\cdot|\det(A_{2})|$.

But the matrix $A_{1}$ is the matrix of evaluation at marked points in $C_{1}$ and $A_{2}$ is the matrix of evaluation at marked points in $C_{2}$. Thus, we have $|\det(A_{1})|=\mult_{\ev_1}(C_1)$ and $|\det(A_{2})|=\mult_{\ev_2}(C_2)$.

Together, we have $\mult_{\ev\times \ft}(C) = \mult_{\ev_{1}}(C_1) \cdot \mult_{\ev_{2}}(C_2) \cdot 2 \cdot (C_{1} \cdot L_{1})_{x_{1}} \cdot (C_{1} \cdot L_{2})_{x_{2}}$.
\end{proof}
We are now ready to prove the main result of this section:
\begin{theorem}[tropical Kontsevich formula for $\F_2$]
\label{kontsevich2}
Let $a$ and $b$ be non-negative integers with $a+b\geq1$. Then
\begin{align*}
\numb{$0$}{$2$}{$a$}{$b$}&=\frac{1}{2}\sum\phi_{2_1}(a_{1},b_{1})\numb{$0$}{$2$}{$a_{1}$}{$b_{1}$}\numb{$0$}{$2$}{$a_{2}$}{$b_{2}$}\\&\quad+\frac{1}{2}\sum \phi_{{2}_{2}}(a_{1},b_{1})\numb{$0$}{$2$}{$a_{1}$}{$b_{1}$}\numb{$0$}{$2$}{$a_{2}$}{$b_{2}$}
\end{align*}
where the first sum goes over all $(a_1,b_1)$ and $(a_2,b_2)$ satisfying \begin{displaymath}(a_1,b_1)+(a_2,b_2)=(a,b),\end{displaymath} $0\leq a_1\leq a$, $0\leq b_1 \leq b$ and $(0,0)\neq (a_1,b_1)\neq (a,b)$, and the second sum goes over all $(a_1,b_1)$ and $(a_2,b_2)$ satisfying \begin{displaymath}(a_1,b_1)+(a_2,b_2)=(a-1,b+2),\end{displaymath} $0\leq a_1\leq a-1$ and $0<b_1<b+2$. We use the shortcuts $\phi_{{2}_{1}}(a_{1},b_{1})$ for
\begin{equation}
\begin{split}
\phi_{{2}_{1}}(a_{1},b_{1})  &= (2a_{1}+b_{1})(2a_{2}+b_{2})(a_{1}b_{2}+a_{2}b_{1}+2a_{1}a_{2})\binom{4a+2b-4}{4a_{1}+2b_{1}-2} \\
& \quad - (2a_{1}+b_{1})(2a_{1}+b_{1})(a_{1}b_{2}+a_{2}b_{1}+2a_{1}a_{2})\binom{4a+2b-4}{4a_{1}+2b_{1}-1}
\end{split}
\label{defphi21}
\end{equation}
and $\phi_{{2}_{2}}(a_{1},b_{1})$ for
\begin{equation}
\begin{split}
\phi_{{2}_{2}}(a_{1},b_{1})  &= 2(2a_{1}+b_{1})(2a_{2}+b_{2})(b_{1}b_{2})\binom{4a+2b-4}{4a_{1}+2b_{1}-2} \\
& \quad - 2(2a_{1}+b_{1})(2a_{1}+b_{1})(b_{1}b_{2})\binom{4a+2b-4}{4a_{1}+2b_{1}-1}.
\end{split}
\label{defphi22}
\end{equation}
\end{theorem}
\begin{proof}
Let $C$ be a curve passing through $L_1$, $L_2$, $p_3,\ldots,p_n$ and mapping to $\lambda_A$ under $\ft$, i.e.\ a curve that contributes to the left hand side of equation \ref{eqdeg} because of remark \ref{rem-generalconditions}. 
By lemma \ref{lem-determinante} and notation \ref{def-mult}, it has to be counted with multiplicity $\mult_{\ev\times \ft}(C)$.
We will show that $C$ can be interpreted as a reducible curve, and that its multiplicity $\mult_{\ev\times \ft}(C)$ can be split into factors corresponding to the irreducible components.

Because of lemma \ref{lem-contractededge} we know that $C$ either has a contracted bounded edge or a string that can be moved to the right as in remark \ref{rem-string}.

If it has a contracted bounded edge, then it is possible that this edge is adjacent to the marked ends $x_1$ and $x_2$. Then the two marked ends are contracted to the same point in $\R^2$, which has to be an intersection point of $L_1$ and $L_2$. Let us call this point $p$. 
Let $C'$ denote the curve that arises after forgetting the marked point $x_1$.
Analogously to 5.5.a) of \cite{GM053} we can show that $\mult_{\ev\times \ft}(C)=\mult_{\ev'}(C')\cdot (L_1.L_2)_p$, where $\ev'$ now denotes the evaluation of $x_2$ at a point combined with all other point evaluations and $(L_1.L_2)_p$ denotes the intersection product of $L_1$ and $L_2$ at $p$.
Instead of counting those curves $C$, we can hence count curves $C'$ meeting the points $p_3,\ldots,p_n$ and an intersection point of $L_1$ and $L_2$. Since $(L_1.L_2)=2$ by the tropical B\'ezout's theorem (4.2 of \cite{RST03}), we can conclude that those curves $C$ contribute $2\numb{$0$}{$2$}{$a$}{$b$}$ to the left hand side of equation \ref{eqdeg}.

If $x_1$ and $x_2$ are not adjacent to the contracted bounded edge, then there have to be bounded edges to both sides of the contracted bounded edge, since all other marked points have to meet different points.
If there are bounded edges on both sides of the contracted bounded edge, we can cut the bounded edge thus producing a reducible curve with two new marked points $z_1$ and $z_2$. Let us call the two components $C_1$ and $C_2$. 
Since $C$ maps to $\lambda_A$ under $\ft$, $x_1$ and $x_2$ have to be on $C_1$.
Let us call the degree of $C_1$ $\Delta_{\F_{2}}(a_1,b_1)$ and the degree of $C_2$ $\Delta_{\F_{2}}(a_2,b_2)$, then we must have $(a_1,b_1)+(a_2,b_2)=(a,b)$, $0\leq a_1\leq a$, $0\leq b_1\leq b$ and $(0,0)\neq (a_1,b_1)\neq (a,b)$.
Analogously to 5.5.b) of \cite{GM053}, we can forget $x_1$ and $x_2$ thus producing a factor of $(C_1.L_1)_{x_1}\cdot (C_1.L_2)_{x_2}$. (By abuse of notation, we use the $x_i$ here for the point in $\R^2$ to which the end $x_i$ is contracted to.)
With the same arguments as in 5.5.b) of \cite{GM053}, we can see that 
\begin{displaymath} \mult_{\ev\times \ft}(C)= \mult_{\ev_1}(C_1) \mult_{\ev_2}(C_2) (C_{1} \cdot C_{2})_{z_{1}=z_{2}} (C_1.L_1)_{x_1}\cdot (C_1.L_2)_{x_2},\end{displaymath} where $\mult_{\ev_1}(C_1)$ denotes the multiplicity of the evaluation at the points on $C_1$. 
Note that $4a_1+2b_1-1$ of the other marked points have to be on $C_1$.
Now instead of counting the curves $C$ with a contracted bounded edge and bounded edges on both sides, we can pick $4a_1+2b_1-1$ of the points $p_5,\ldots,p_n$ ($\binom{4a+2b-4}{4a_1+2b_1-1}$ possibilities) and count curves $C_1$ through those points, and $C_2$ through the remaining points. Again by tropical B\'ezout's theorem we have $(C_1.L_1)=(2a_1+b_1)$ choices to attach $x_1$ and $(C_1.L_2)=(2a_1+b_1)$ choices to attach $x_2$, and we have $(C_1.C_2)=(a_1b_2+a_2b_1+2a_1a_2)$ choices to glue $C_1$ and $C_2$ to a possible $C$. Thus those curves contribute 
\begin{displaymath}\sum(2a_{1}+b_{1})(2a_{1}+b_{1})(a_{1}b_{2}+a_{2}b_{1}+2a_{1}a_{2})\binom{4a+2b-4}{4a_{1}+2b_{1}-1}\numb{$0$}{$2$}{$a_{1}$}{$b_{1}$}\numb{$0$}{$2$}{$a_{2}$}{$b_{2}$},\end{displaymath}
where the sum goes over all $(a_1,b_1)+(a_2,b_2)=(a,b)$, $0\leq a_1\leq a$, $0\leq b_1\leq b$ and $(0,0)\neq (a_1,b_1)\neq (a,b)$. In the formula we want to prove, we can see this contribution negatively on the right hand side.

Finally, if $C$ has a string as in remark \ref{rem-string}, then by lemma \ref{lem-multstring} we can conclude that $C$ contributes
\[ \mult_{\ev\times \ft}(C) = \mult_{\ev_1}(C_1) \cdot \mult_{\ev_2}(C_2) \cdot 2 \cdot (C_{1} \cdot L_{1})_{x_{1}} \cdot (C_{1} \cdot L_{2})_{x_{2}}. \]
Instead of counting such curves $C$, we can pick $4a_1+2b_1-1$ of the points $p_5,\ldots,p_n$ ($\binom{4a+2b-4}{4a_1+2b_1-1}$ possibilities) and count curves $C_1$ of degree $\Delta_{\F_{2}}(a_1,b_1)$ through those points, and $C_2$ of degree $\Delta_{\F_{2}}(a_2,b_2)$ through the remaining points, where now $(a_1,b_1)+(a_2,b_2)=(a-1,b+2)$. 
There are again $(2a_1+b_1)$ possibilities to attach $x_1$ and also $(2a_1+b_1)$ possibilities to attach $x_2$ to $C_1$. There are $b_1b_2$ choices to pick the edges $E'$ and $E''$ to which we can attach the string $S$.
Hence those curves contribute
\begin{displaymath}\sum 2(2a_{1}+b_{1})(2a_{1}+b_{1})(b_1b_2)\binom{4a+2b-4}{4a_{1}+2b_{1}-1}\numb{$0$}{$2$}{$a_{1}$}{$b_{1}$}\numb{$0$}{$2$}{$a_{2}$}{$b_{2}$},\end{displaymath}
where now the sum goes over all $(a_1,b_1)+(a_2,b_2)=(a-1,b+2)$, $0\leq a_1\leq a-1$, $0< b_1< b+2$ and $(0,0)\neq (a_1,b_1)\neq (a,b)$. In the formula we want to prove, this contribution appears negatively on the right hand side.

Performing the same analysis for the right hand side of equation \ref{eqdeg} and collecting the terms to the different sides, the statement follows. 

\end{proof}
The following formula for $\F_0$ can be proved analogously. The proof is easier in fact, since all curves have a contracted bounded edge and the special case of curves having a string that can be moved to the right as in remark \ref{rem-string} does not occur here. We therefore skip the proof.
For more details, see \cite{Fra08}.

\begin{theorem}[tropical Kontsevich formula for $\F_0$]
\label{kontsevich0}
Let $a$ and $b$ be non-negative integers with $a+b\geq1$. Then
\begin{equation*}
\numb{$0$}{$0$}{$a$}{$a+b$}=\frac{1}{2}\sum\phi_{0}(a_{1},a_{1}+b_{1})\numb{$0$}{$0$}{$a_{1}$}{$a_{1}+b_{1}$}\numb{$0$}{$0$}{$a_{2}$}{$a_{2}+b_{2}$}
\end{equation*}
where the sum goes over all $(a_1,b_1)$ and $(a_2,b_2)$ satisfying \begin{displaymath}(a_{1},a_{1}+b_{1})+(a_{2},a_{2}+b_{2})=(a,a+b),\end{displaymath} $0\leq a_{1} \leq a$, $-a_{1} \leq b_{1} \leq b+a_{2}$ and $(0,0)\neq(a_{1},a_{1}+b_{1})\neq(a,a+b)$. We use the shortcut $\phi_{0}(a_{1},a_{1}+b_{1})$ for
\begin{equation}
\begin{split}
\label{defphi0}
\phi_{0}(a_{1}&,a_{1}+b_{1})  = (2a_{1}+b_{1})(2a_{2}+b_{2})(a_{1}b_{2}+a_{2}b_{1}+2a_{1}a_{2})\binom{4a+2b-4}{4a_{1}+2b_{1}-2} \\
& \quad - (2a_{1}+b_{1})(2a_{1}+b_{1})(a_{1}b_{2}+a_{2}b_{1}+2a_{1}a_{2})\binom{4a+2b-4}{4a_{1}+2b_{1}-1}.
\end{split}
\end{equation}
\end{theorem}

\section{The proof of theorem \ref{thm-main}}
\label{sec-2}


Recall the equation from theorem \ref{thm-main} we want to prove.
Note that for $k=a$ we would obtain a summand $\binom{b+2a}{a}\mathcal{N}^{0}_{\F_{2}}(0,(b+2a))$ which is $0$ for all $a$, $b\in\Z_{\geq0}$ with $a+b\geq1$ except for $a=0$ and $b=1$. As in this special case the statement still holds, we may add this summand for $k=a$ and deal with the slightly modified equation
\begin{displaymath}
\numb{$0$}{$0$}{$a$}{$a+b$} = \sum_{k=0}^{a}{\binom{b+2k}{k}\numb{$0$}{$2$}{$a-k$}{$b+2k$}}.
\end{displaymath}
We will see later that this is useful.

\begin{example}
Let us consider the formula for small $a$ in more detail. For $a=0$ and $a=1$ we have $\numb{$0$}{$0$}{$a$}{$a+b$}=\numb{$0$}{$2$}{$a$}{$b$}$ for all $b\in\Z_{\geq 0}$. Hence, in these cases the Gromov-Witten invariants are enumerative for $\F_{2}$ as well. For $a=2$ and $b=0$ we have $\numb{$0$}{$2$}{$2$}{$0$}=10$ while the associated Gromov-Witten invariant is $\numb{$0$}{$0$}{$2$}{$2$}=12$. This is the first interesting case. The formula gives an interpretation of the difference as multiple of an enumerative number of curves of different class:
\[ \numb{$0$}{$0$}{$2$}{$2$} = 12 = 10 + 2\times 1 = \numb{$0$}{$2$}{$2$}{$0$} + \binom{2}{1}\numb{$0$}{$2$}{$1$}{$2$} + \underbrace{\binom{4}{2}\numb{$0$}{$2$}{$0$}{$4$}}_{=0} \text{.} \]
\end{example}

We need the following combinatorial indentity involving binomial coefficients for our proof.

\begin{lemma}
\label{binomgleichung}
Let $n$, $m$, $k\in\N$. Then
\[\sum_{i=0}^{k}(mi+n(k-i)-2i(k-i))\binom{n}{i}\binom{m}{k-i}=2\cdot n\cdot m \cdot \binom{n+m-2}{k-1}.\]
\end{lemma}

\begin{proof}

The equation is essentially a consequence of Vandermonde's identity which states that
\[\sum_{i=0}^{k}\binom{n}{i}\binom{m}{k-i}=\binom{n+m}{k}.\]
Using this we have
\[\sum_{i=0}^{k}mi\binom{n}{i}\binom{m}{k-i} = \sum_{i=0}^{k}nm\binom{n-1}{i-1}\binom{m}{k-i} = nm\binom{n+m-1}{k-1}\] and
\[\sum_{i=0}^{k}n(k-i)\binom{n}{i}\binom{m}{k-i} = \sum_{i=0}^{k}nm\binom{n}{i}\binom{m-1}{k-i-1} = nm\binom{n+m-1}{k-1}\]
and
\begin{equation*}
\begin{split}
\sum_{i=0}^{k}(-2i(k-i))\binom{n}{i}\binom{m}{k-i} &= -2\sum_{i=0}^{k}nm\binom{n-1}{i-1}\binom{m-1}{k-i-1} \\
&= -2nm\binom{n+m-2}{k-2}
\end{split}
\end{equation*}
as $i\binom{n}{i}=n\binom{n-1}{i-1}$ and $(k-i)\binom{m}{k-i}=m\binom{m-1}{k-i-1}$ and thus
\begin{equation*}
\begin{split}
& \sum_{i=0}^{k}(mi+n(k-i)-2i(k-i))\binom{n}{i}\binom{m}{k-i} \\
& \quad = 2nm\left( \binom{n+m-1}{k-1}-\binom{n+m-2}{k-2} \right) \\
& \quad = 2nm\binom{n+m-2}{k-1}
\end{split}
\end{equation*}
where the last equality follows by Pascal's rule.

\end{proof}

\begin{proof}[Proof of theorem \ref{thm-main}:]
We will prove that the statement holds for all integers $a\geq0$ and $b\geq-a$ with $2a+b\geq1$ by induction on $2a+b$. Note that $b$ may be negative. But as $b \geq -a$ and $2a+b \geq 1$ we have $a+b \geq 0$ and hence the left hand side is well defined. On the right hand side we may have negative entries. We define $\numb{$0$}{$2$}{$a$}{$b$}$ to be $0$ for all $b<0$, $a\geq 0$. In particular, we conclude that the statement holds for all non-negative integers $a$ and $b$ with $a+b\geq1$.

The induction beginning for $a=0$ and $b=1$ resp.\ $a=1$ and $b=-1$ is straight forward, we need to use the extra summand with $k=a$ however.


Now let $a\geq0$ and $b\geq-a$ be integers such that $2a+b \geq 1$.

We can assume that the statement holds for all integers $a_{i}\geq0$ and $b_{i}\geq-a_{i}$ with $1 \leq 2a_{i}+b_{i}<2a+b$.

First let us consider the left hand side. We know by theorem \ref{kontsevich0} that
\begin{equation*}
 2\numb{$0$}{$0$}{$a$}{$a+b$}  =
\sum \phi_{0}(a_{1},a_{1}+b_{1})\numb{$0$}{$0$}{$a_{1}$}{$a_{1}+b_{1}$}\numb{$0$}{$0$}{$a_{2}$}{$a_{2}+b_{2}$}
\end{equation*}
 where the sum goes over all $(a_1,b_1)$ and $(a_2,b_2)$ satisfying $(a_{1},a_{1}+b_{1})+(a_{2},a_{2}+b_{2})=(a,a+b)$, $0\leq a_{1} \leq a$, $-a_{1} \leq b_{1} \leq b+a_{2}$ and $(0,0)\neq(a_{1},a_{1}+b_{1})\neq(a,a+b)$, and $\phi_{0}(a_{1},a_{1}+b_{1})$ is defined by equation \ref{defphi0}.

As $2a_{1}+b_{1}<2a+b$ and $2a_{2}+b_{2}<2a+b$ for all $a_{1}$, $a_{2}$, $b_{1}$ and $b_{2}$ we have by the induction hypothesis that
$\numb{$0$}{$0$}{$a_{1}$}{$a_{1}+b_{1}$} = \sum_{i=0}^{a_{1}}{\binom{b_{1}+2i}{i}\numb{$0$}{$2$}{$a_{1}-i$}{$b_{1}+2i$}}$ and
$\numb{$0$}{$0$}{$a_{2}$}{$a_{2}+b_{2}$} = \sum_{j=0}^{a_{2}}{\binom{b_{2}+2j}{j}\numb{$0$}{$2$}{$a_{2}-j$}{$b_{2}+2j$}}$.
Hence we have
\begin{equation*}
\begin{split}
 2\numb{$0$}{$0$}{$a$}{$a+b$}  &=
\sum \phi_{0}(a_{1},a_{1}+b_{1}) \cdot \left(\sum_{i=0}^{a_{1}}{\binom{b_{1}+2i}{i}\numb{$0$}{$2$}{$a_{1}-i$}{$b_{1}+2i$}}\right) \\
& \qquad \hspace{5em} \cdot \left(\sum_{j=0}^{a_{2}}{\binom{b_{2}+2j}{j}\numb{$0$}{$2$}{$a_{2}-j$}{$b_{2}+2j$}}\right) \\
&=
\sum_{k=0}^{a}
\sum_{\substack{i=0 \\ j=k-i}}^{k}  \sum
\bigg(\binom{b_{1}+2i}{i}\binom{b_{2}+2j}{j}
\phi_{0}(a_{1},a_{1}+b_{1}) \\
& \qquad \hspace{5em}  \cdot \numb{$0$}{$2$}{$a_{1}-i$}{$b_{1}+2i$}\numb{$0$}{$2$}{$a_{2}-j$}{$b_{2}+2j$}\bigg)
\end{split}
\end{equation*}

Let us consider the range of $a_{1}$ and $b_{1}$ in the third sum for a fixed $k$ and $i$. As $\numb{$0$}{$2$}{$a_{1}-i$}{$b_{1}+2i$}=0$ for all $0 \leq a_{1} < i$ and $\numb{$0$}{$2$}{$a_{2}-j$}{$b_{2}+2j$}=0$ for all $0 \leq a_{2} < j$ (i.e.\ for all $a-0 = a \geq a-a_{2} = a_{1} > a-j = a-k+i$) we can forget about the summands where $0 \leq a_{1} < i$ or $a-k+i < a_{1} \leq a$. As $\binom{b_{1}+2i}{i}=0$ for all $b_{1} < -i$ and $\binom{b_{2}+2j}{j}=0$ for all $b_{2} < -j$ (i.e.\ for all $b- b_{2} = b_{1} > b+j$) the range of those $b_{1}$ which give a contribution is $-i \leq b_{1} \leq b+j$. We may add summands for $-2i \leq b_{1} < -i$ and $b+j < b_{1} \leq b+2j$ since they are $0$ anyway. Hence we may restrict our attention to those $(a_{1},a_{1}+b_{1})+(a_{2},a_{2}+b_{2})=(a,a+b)$ with
$i \leq a_{1} \leq a-k+i$ and $ -2i \leq b_{1} \leq b+2j$ such that $ (0,0)\neq(a_{1},a_{1}+b_{1})\neq(a,a+b)$.
With the definitions
\begin{displaymath}
\begin{array}{ll}
a_{1}':=a_{1}-i &  \quad a_{2}':= (a-k)-a_{1}'= a_{2}-j \\
b_{1}':=b_{1}+2i & \quad b_{2}':= (b+2k)-b_{1}'= b_{2}+2j
\end{array}
\end{displaymath}
this is equivalent to considering all pairs \[(a_{1}',b_{1}')+(a_{2}',b_{2}')=(a-k,b+2k)\] with \[0 \leq a_{1}' \leq a-k \text{ and } 0 \leq b_{1}' \leq b+2k \text{ such that } (0,0)\neq(a_{1}',b_{1}')\neq(a,b) \text{.}\]
Let the sums in the following equation go over those pairs now. Then


\begin{equation*}
\begin{split}
2\numb{$0$}{$0$}{$a$}{$a+b$} &= \sum_{k=0}^{a}\sum_{\substack{i=0 \\ j=k-i}}^{k}\sum
\bigg( \binom{b_{1}'}{i}\binom{b_{2}'}{j}
\phi_{0}(a_{1}'+i,a_{1}'+b_{1}'-i) \\
&  \quad \hspace{5em} \cdot \numb{$0$}{$2$}{$a_{1}'$}{$b_{1}'$}\numb{$0$}{$2$}{$a_{2}'$}{$b_{2}'$}\bigg) \\
& = \sum_{k=0}^{a}\sum
\bigg(\sum_{\substack{i=0 \\ j=k-i}}^{k}
\binom{b_{1}'}{i}\binom{b_{2}'}{j}
\phi_{0}(a_{1}'+i,a_{1}'+b_{1}'-i) \\
& \quad \hspace{5em} \cdot \numb{$0$}{$2$}{$a_{1}'$}{$b_{1}'$}\numb{$0$}{$2$}{$a_{2}'$}{$b_{2}'$}\bigg)
\end{split}
\end{equation*}
where
\begin{equation*}
\begin{split}
& \phi_{0}(a_{1}'+i,a_{1}'+b_{1}'-i) \\
& \quad = (2a_{1}'+b_{1}')(2a_{2}'+b_{2}')(a_{1}'b_{2}'+a_{2}'b_{1}'+2a_{1}'a_{2}'\\
& \qquad \hspace{10em} +b_{2}'i+b_{1}'j-2ij)\binom{4a+2b-4}{4a_{1}'+2b_{1}'-2} \\
& \qquad - (2a_{1}'+b_{1}')(2a_{1}'+b_{1}')(a_{1}'b_{2}'+a_{2}'b_{1}'+2a_{1}'a_{2}'\\
& \qquad \hspace{10em} +b_{2}'i+b_{1}'j-2ij)\binom{4a+2b-4}{4a_{1}'+2b_{1}'-1} \\
& \quad =(2a_{1}'+b_{1}')(2a_{2}'+b_{2}')(a_{1}'b_{2}'+a_{2}'b_{1}'+2a_{1}'a_{2}')
\binom{4a+2b-4}{4a_{1}'+2b_{1}'-2} \\
& \qquad \hspace{5em} - (2a_{1}'+b_{1}')(2a_{1}'+b_{1}')(a_{1}'b_{2}'+a_{2}'b_{1}'+2a_{1}'a_{2}')
\binom{4a+2b-4}{4a_{1}'+2b_{1}'-1} \\
& \qquad + (2a_{1}'+b_{1}')(2a_{2}'+b_{2}')(b_{2}'i+b_{1}'j-2ij)
\binom{4a+2b-4}{4a_{1}'+2b_{1}'-2} \\
& \qquad \hspace{5em} - (2a_{1}'+b_{1}')(2a_{1}'+b_{1}')(b_{2}'i+b_{1}'j-2ij)
\binom{4a+2b-4}{4a_{1}'+2b_{1}'-1}
\end{split}
\end{equation*}
by equation (\ref{defphi0}) in theorem \ref{kontsevich0}.

Let us stop here and consider the right hand side of the equation we want to prove. We have by Theorem \ref{kontsevich2}
\begin{equation*}
\begin{split}
& \sum_{k=0}^{a}{\binom{b+2k}{k}2\numb{$0$}{$2$}{$a-k$}{$b+2k$}}  \\ =& \sum_{k=0}^{a} \binom{b+2k}{k}
\cdot\left(\sum{\phi_{{2}_{1}}(a_{1},b_{1}) \numb{$0$}{$2$}{$a_{1}$}{$b_{1}$}\numb{$0$}{$2$}{$a_{2}$}{$b_{2}$}}\right. \\
& \quad \hspace{5em} + \left. \sum{\phi_{{2}_{2}}(a_{1},b_{1})\numb{$0$}{$2$}{$a_{1}$}{$b_{1}$}\numb{$0$}{$2$}{$a_{2}$}{$b_{2}$}}\right) \\
\end{split}
\end{equation*}
where the first sum goes over all pairs such that $(a_1,b_1)+(a_2,b_2)=(a-k,b+2k)$ and the second sum goes over all pairs such that $(a_1,b_1)+(a_2,b_2)=(a-(k+1),b+2(k+1))$. We use the shortcuts $\phi_{{2}_{1}}(a_{1},b_{1})$ and $\phi_{{2}_{2}}(a_{1},b_{1})$ as defined in equation \ref{defphi21} and
\ref{defphi22}.

Since for $k=0$ the binomial coefficient $\binom{b+2(k-1)}{k-1}$ is $0$ and for $k=a$ there are no $a_{1}$ and $b_{1}$ which satisfy $(a_1,b_1)+(a_2,b_2)=(a-(k+1),b+2(k+1))$ we can merge the two sums and get
\begin{equation*}
\begin{split}
\sum_{k=0}^{a}
\sum & \bigg(\bigg( \binom{b+2k}{k}\phi_{{2}_{1}}(a_{1},b_{1}) +\binom{b+2(k-1)}{k-1}\phi_{{2}_{2}}(a_{1},b_{1})\bigg) \\
& \quad \cdot  \numb{$0$}{$2$}{$a_{1}$}{$b_{1}$}\numb{$0$}{$2$}{$a_{2}$}{$b_{2}$}\bigg)
\end{split}
\end{equation*}
where the sum now goes over all pairs such that $(a_1,b_1)+(a_2,b_2)=(a-k,b+2k)$.

Thus it remains to show that
\begin{equation*}
\begin{split}
&\sum_{k=0}^{a}\sum\bigg(\sum_{\substack{i=0 \\ j=k-i}}^{k}
\binom{b_{1}}{i}
\binom{b_{2}}{j}
\phi_{0}(a_{1}+i,a_{1}+b_{1}-i)\bigg)  \cdot \numb{$0$}{$2$}{$a_{1}$}{$b_{1}$}\numb{$0$}{$2$}{$a_{2}$}{$b_{2}$} \\
 =&
\sum_{k=0}^{a}
\sum \bigg( \binom{b+2k}{k}\phi_{{2}_{1}}(a_{1},b_{1})+\binom{b+2(k-1)}{k-1}\phi_{{2}_{2}}(a_{1},b_{1})\bigg) \\
&
\hspace{8em} \cdot \numb{$0$}{$2$}{$a_{1}$}{$b_{1}$}\numb{$0$}{$2$}{$a_{2}$}{$b_{2}$}.
\end{split}
\end{equation*}

Therefore we will show that
\begin{equation*}
\begin{split}
& \sum_{\substack{i=0 \\ j=k-i}}^{k}\binom{b_{1}}{i}\binom{b_{2}}{j}\phi_{0}(a_{1}+i,a_{1}+b_{1}-i) \\
=&
\binom{b+2k}{k}\phi_{{2}_{1}}(a_{1},b_{1})+\binom{b+2(k-1)}{k-1}\phi_{{2}_{2}}(a_{1},b_{1})
\end{split}
\label{toshow2}
\end{equation*}
for all $k\in\{0,...,a\}$ and for all integers $0\leq a_{1}$, $a_{2}\leq a-k$, $0 \leq b_{1}$, $b_{2} \leq b+2k$ with $a_{1}+a_{2}=a-k$, $b_{1}+b_{2}=b+2k$ and $(0,0)\neq(a_{1},b_{1})\neq(a-k,b+2k)$.

We use the identity from Lemma \ref{binomgleichung}.

It is
\begin{equation*}
\begin{split}
& \sum_{\substack{i=0 \\ j=k-i}}^{k}\binom{b_{1}}{i}\binom{b_{2}}{j}\phi_{0}(a_{1}+i,a_{1}+b_{1}-i) \\
\stackrel{(\ref{defphi0})}{=} &
\sum_{\substack{i=0 \\ j=k-i}}^{k}\bigg[\binom{b_{1}}{i}\binom{b_{2}}{j} \\
& \quad \hspace{5em} \cdot
\bigg(
(2a_{1}+b_{1})(2a_{2}+b_{2})(a_{1}b_{2}+a_{2}b_{1}+2a_{1}a_{2})
\binom{4a+2b-4}{4a_{1}+2b_{1}-2} \\
& \qquad \hspace{5em} - (2a_{1}+b_{1})(2a_{1}+b_{1})(a_{1}b_{2}+a_{2}b_{1}+2a_{1}a_{2})
\binom{4a+2b-4}{4a_{1}+2b_{1}-1} \\
& \qquad \hspace{5em} + (2a_{1}+b_{1})(2a_{2}+b_{2})(b_{2}i+b_{1}j-2ij)
\binom{4a+2b-4}{4a_{1}+2b_{1}-2} \\
& \qquad \hspace{5em} - (2a_{1}+b_{1})(2a_{1}+b_{1})(b_{2}i+b_{1}j-2ij)
\binom{4a+2b-4}{4a_{1}+2b_{1}-1}
\bigg)\bigg]
\end{split}
\end{equation*}
\begin{equation*}
\begin{split}
& \quad =
\bigg(\sum_{\substack{i=0 \\ j=k-i}}^{k}\binom{b_{1}}{i}\binom{b_{2}}{j}\bigg) \\
& \quad \hspace{5em} \cdot \bigg(
(2a_{1}+b_{1})(2a_{2}+b_{2})(a_{1}b_{2}+a_{2}b_{1}+2a_{1}a_{2})
\binom{4a+2b-4}{4a_{1}+2b_{1}-2} \\
& \qquad \hspace{5em} - (2a_{1}+b_{1})(2a_{1}+b_{1})(a_{1}b_{2}+a_{2}b_{1}+2a_{1}a_{2})
\binom{4a+2b-4}{4a_{1}+2b_{1}-1}\bigg) \\
& \qquad + \bigg(\sum_{\substack{i=0 \\ j=k-i}}^{k}\binom{b_{1}}{i}\binom{b_{2}}{j}(b_{2}i+b_{1}j-2ij)\bigg) \\
& \quad \hspace{5em} \cdot \bigg((2a_{1}+b_{1})(2a_{2}+b_{2})
\binom{4a+2b-4}{4a_{1}+2b_{1}-2} \\
& \qquad \hspace{5em} - (2a_{1}+b_{1})(2a_{1}+b_{1})
\binom{4a+2b-4}{4a_{1}+2b_{1}-1}
\bigg)
\end{split}
\end{equation*}
\begin{equation*}
\begin{split}
& \quad = \binom{b_{1}+b_{2}}{i+j}
\bigg(
(2a_{1}+b_{1})(2a_{2}+b_{2})(a_{1}b_{2}+a_{2}b_{1}+2a_{1}a_{2})
\binom{4a+2b-4}{4a_{1}+2b_{1}-2} \\
& \qquad \hspace{5em} - (2a_{1}+b_{1})(2a_{1}+b_{1})(a_{1}b_{2}+a_{2}b_{1}+2a_{1}a_{2})
\binom{4a+2b-4}{4a_{1}+2b_{1}-1}\bigg) \\
& \qquad + \bigg(2b_{1}b_{2}\binom{b_{1}+b_{2}-2}{i+j-1}\bigg)
\bigg((2a_{1}+b_{1})(2a_{2}+b_{2})
\binom{4a+2b-4}{4a_{1}+2b_{1}-2} \\
& \quad \hspace{12em} - (2a_{1}+b_{1})(2a_{1}+b_{1})
\binom{4a+2b-4}{4a_{1}+2b_{1}-1}\bigg)
\end{split}
\end{equation*}
\begin{equation*}
\begin{split}
& \quad =
\binom{b_{1}+b_{2}}{i+j}
\bigg((2a_{1}+b_{1})(2a_{2}+b_{2})(a_{1}b_{2}+a_{2}b_{1}+2a_{1}a_{2})
\binom{4a+2b-4}{4a_{1}+2b_{1}-2} \\
& \qquad \hspace{5em} - (2a_{1}+b_{1})(2a_{1}+b_{1})(a_{1}b_{2}+a_{2}b_{1}+2a_{1}a_{2})
\binom{4a+2b-4}{4a_{1}+2b_{1}-1}\bigg) \\
& \qquad + \binom{b_{1}+b_{2}-2}{i+j-1}
\bigg(2(2a_{1}+b_{1})(2a_{2}+b_{2})b_{1}b_{2}
\binom{4a+2b-4}{4a_{1}+2b_{1}-2} \\
& \quad \hspace{12em} - 2(2a_{1}+b_{1})(2a_{1}+b_{1})b_{1}b_{2}
\binom{4a+2b-4}{4a_{1}+2b_{1}-1}\bigg) \\
& \quad =
\binom{b+2k}{k}\phi_{{2}_{1}}(a_{1},b_{1})+\binom{b+2(k-1)}{k-1}\phi_{{2}_{2}}(a_{1},b_{1})
\end{split}
\end{equation*}
where the last equality follows from equation (\ref{defphi21}) and (\ref{defphi22}).
This completes the proof.
\end{proof}


\end{document}

%% file: polytope.pstex_t
\begin{picture}(0,0)%
\includegraphics{polytope.pstex}%
\end{picture}%
\setlength{\unitlength}{3191sp}%
\begingroup\makeatletter\ifx\SetFigFont\undefined%
\gdef\SetFigFont#1#2#3#4#5{%
  \reset@font\fontsize{#1}{#2pt}%
  \fontfamily{#3}\fontseries{#4}\fontshape{#5}%
  \selectfont}%
\fi\endgroup%
\begin{picture}(3110,2032)(2611,-2981)
\put(5491,-2491){\makebox(0,0)[lb]{\smash{{\SetFigFont{8}{9.6}{\familydefault}{\mddefault}{\updefault}{\color[rgb]{0,0,0}$b$}%
}}}}
\put(4951,-2941){\makebox(0,0)[lb]{\smash{{\SetFigFont{8}{9.6}{\familydefault}{\mddefault}{\updefault}{\color[rgb]{0,0,0}$a$}%
}}}}
\put(3376,-2941){\makebox(0,0)[lb]{\smash{{\SetFigFont{8}{9.6}{\familydefault}{\mddefault}{\updefault}{\color[rgb]{0,0,0}$a$}%
}}}}
\put(5176,-1636){\makebox(0,0)[lb]{\smash{{\SetFigFont{8}{9.6}{\familydefault}{\mddefault}{\updefault}{\color[rgb]{0,0,0}$a$}%
}}}}
\put(2611,-2266){\makebox(0,0)[lb]{\smash{{\SetFigFont{8}{9.6}{\familydefault}{\mddefault}{\updefault}{\color[rgb]{0,0,0}$a+b$}%
}}}}
\put(4096,-2221){\makebox(0,0)[lb]{\smash{{\SetFigFont{8}{9.6}{\familydefault}{\mddefault}{\updefault}{\color[rgb]{0,0,0}$2a+b$}%
}}}}
\end{picture}%

%% file: curvebsp.pstex_t
\begin{picture}(0,0)%
\includegraphics{curvebsp.pstex}%
\end{picture}%
\setlength{\unitlength}{3947sp}%
\begingroup\makeatletter\ifx\SetFigFont\undefined%
\gdef\SetFigFont#1#2#3#4#5{%
  \reset@font\fontsize{#1}{#2pt}%
  \fontfamily{#3}\fontseries{#4}\fontshape{#5}%
  \selectfont}%
\fi\endgroup%
\begin{picture}(4008,1895)(5326,-5961)
\put(6784,-4744){\makebox(0,0)[lb]{\smash{{\SetFigFont{7}{8.4}{\familydefault}{\mddefault}{\updefault}{\color[rgb]{0,0,0}$h$}%
}}}}
\put(5542,-5932){\makebox(0,0)[lb]{\smash{{\SetFigFont{7}{8.4}{\familydefault}{\mddefault}{\updefault}{\color[rgb]{0,0,0}$\Gamma$}%
}}}}
\put(9052,-4204){\makebox(0,0)[lb]{\smash{{\SetFigFont{7}{8.4}{\familydefault}{\mddefault}{\updefault}{\color[rgb]{0,0,0}$\R^2$}%
}}}}
\put(5596,-4150){\makebox(0,0)[lb]{\smash{{\SetFigFont{7}{8.4}{\familydefault}{\mddefault}{\updefault}{\color[rgb]{0,0,0}$x_1$}%
}}}}
\put(6190,-4798){\makebox(0,0)[lb]{\smash{{\SetFigFont{7}{8.4}{\familydefault}{\mddefault}{\updefault}{\color[rgb]{0,0,0}$x_2$}%
}}}}
\put(6190,-5446){\makebox(0,0)[lb]{\smash{{\SetFigFont{7}{8.4}{\familydefault}{\mddefault}{\updefault}{\color[rgb]{0,0,0}$x_4$}%
}}}}
\put(5650,-4906){\makebox(0,0)[lb]{\smash{{\SetFigFont{7}{8.4}{\familydefault}{\mddefault}{\updefault}{\color[rgb]{0,0,0}$x_3$}%
}}}}
\put(5326,-5068){\makebox(0,0)[lb]{\smash{{\SetFigFont{7}{8.4}{\familydefault}{\mddefault}{\updefault}{\color[rgb]{0,0,0}$x_5$}%
}}}}
\end{picture}%

%% file: curvebsp2.pstex_t
\begin{picture}(0,0)%
\includegraphics{curvebsp2.pstex}%
\end{picture}%
\setlength{\unitlength}{3947sp}%
\begingroup\makeatletter\ifx\SetFigFont\undefined%
\gdef\SetFigFont#1#2#3#4#5{%
  \reset@font\fontsize{#1}{#2pt}%
  \fontfamily{#3}\fontseries{#4}\fontshape{#5}%
  \selectfont}%
\fi\endgroup%
\begin{picture}(4213,1540)(5326,-5567)
\put(5596,-4150){\makebox(0,0)[lb]{\smash{{\SetFigFont{7}{8.4}{\familydefault}{\mddefault}{\updefault}{\color[rgb]{0,0,0}$x_1$}%
}}}}
\put(6190,-4798){\makebox(0,0)[lb]{\smash{{\SetFigFont{7}{8.4}{\familydefault}{\mddefault}{\updefault}{\color[rgb]{0,0,0}$x_2$}%
}}}}
\put(6190,-5446){\makebox(0,0)[lb]{\smash{{\SetFigFont{7}{8.4}{\familydefault}{\mddefault}{\updefault}{\color[rgb]{0,0,0}$x_4$}%
}}}}
\put(5650,-4906){\makebox(0,0)[lb]{\smash{{\SetFigFont{7}{8.4}{\familydefault}{\mddefault}{\updefault}{\color[rgb]{0,0,0}$x_3$}%
}}}}
\put(5326,-5068){\makebox(0,0)[lb]{\smash{{\SetFigFont{7}{8.4}{\familydefault}{\mddefault}{\updefault}{\color[rgb]{0,0,0}$x_5$}%
}}}}
\put(8551,-5536){\makebox(0,0)[lb]{\smash{{\SetFigFont{7}{8.4}{\familydefault}{\mddefault}{\updefault}{\color[rgb]{0,0,0}$x_3$}%
}}}}
\put(9076,-5536){\makebox(0,0)[lb]{\smash{{\SetFigFont{7}{8.4}{\familydefault}{\mddefault}{\updefault}{\color[rgb]{0,0,0}$x_4$}%
}}}}
\put(9001,-4786){\makebox(0,0)[lb]{\smash{{\SetFigFont{7}{8.4}{\familydefault}{\mddefault}{\updefault}{\color[rgb]{0,0,0}$l_1+l_2$}%
}}}}
\put(8626,-4111){\makebox(0,0)[lb]{\smash{{\SetFigFont{7}{8.4}{\familydefault}{\mddefault}{\updefault}{\color[rgb]{0,0,0}$x_1$}%
}}}}
\put(9001,-4111){\makebox(0,0)[lb]{\smash{{\SetFigFont{7}{8.4}{\familydefault}{\mddefault}{\updefault}{\color[rgb]{0,0,0}$x_2$}%
}}}}
\put(7246,-4150){\makebox(0,0)[lb]{\smash{{\SetFigFont{7}{8.4}{\familydefault}{\mddefault}{\updefault}{\color[rgb]{0,0,0}$x_1$}%
}}}}
\put(7840,-4798){\makebox(0,0)[lb]{\smash{{\SetFigFont{7}{8.4}{\familydefault}{\mddefault}{\updefault}{\color[rgb]{0,0,0}$x_2$}%
}}}}
\put(7840,-5446){\makebox(0,0)[lb]{\smash{{\SetFigFont{7}{8.4}{\familydefault}{\mddefault}{\updefault}{\color[rgb]{0,0,0}$x_4$}%
}}}}
\put(7060,-4892){\makebox(0,0)[lb]{\smash{{\SetFigFont{7}{8.4}{\familydefault}{\mddefault}{\updefault}{\color[rgb]{0,0,0}$x_3$}%
}}}}
\put(7610,-5059){\makebox(0,0)[lb]{\smash{{\SetFigFont{7}{8.4}{\familydefault}{\mddefault}{\updefault}{\color[rgb]{0,0,0}$l_1$}%
}}}}
\put(7512,-5237){\makebox(0,0)[lb]{\smash{{\SetFigFont{7}{8.4}{\familydefault}{\mddefault}{\updefault}{\color[rgb]{0,0,0}$l_2$}%
}}}}
\end{picture}%

%% file: detbspcomb.pstex_t
\begin{picture}(0,0)%
\includegraphics{detbspcomb.pstex}%
\end{picture}%
\setlength{\unitlength}{3947sp}%
\begingroup\makeatletter\ifx\SetFigFont\undefined%
\gdef\SetFigFont#1#2#3#4#5{%
  \reset@font\fontsize{#1}{#2pt}%
  \fontfamily{#3}\fontseries{#4}\fontshape{#5}%
  \selectfont}%
\fi\endgroup%
\begin{picture}(1865,1687)(3826,-4389)
\put(4465,-3261){\makebox(0,0)[lb]{\smash{{\SetFigFont{9}{10.8}{\familydefault}{\mddefault}{\updefault}{\color[rgb]{0,0,0}$\binom{1}{0}$}%
}}}}
\put(5295,-3453){\makebox(0,0)[lb]{\smash{{\SetFigFont{9}{10.8}{\familydefault}{\mddefault}{\updefault}{\color[rgb]{0,0,0}$\binom{2}{1}$}%
}}}}
\put(4912,-3644){\makebox(0,0)[lb]{\smash{{\SetFigFont{9}{10.8}{\familydefault}{\mddefault}{\updefault}{\color[rgb]{0,0,0}$\binom{-1}{-1}$}%
}}}}
\put(4784,-3964){\makebox(0,0)[lb]{\smash{{\SetFigFont{9}{10.8}{\familydefault}{\mddefault}{\updefault}{\color[rgb]{0,0,0}$\binom{0}{-1}$}%
}}}}
\put(4209,-4092){\makebox(0,0)[lb]{\smash{{\SetFigFont{9}{10.8}{\familydefault}{\mddefault}{\updefault}{\color[rgb]{0,0,0}$\binom{-1}{0}$}%
}}}}
\put(3826,-4155){\makebox(0,0)[lb]{\smash{{\SetFigFont{9}{10.8}{\familydefault}{\mddefault}{\updefault}{\color[rgb]{0,0,0}$x_4$}%
}}}}
\put(5232,-4347){\makebox(0,0)[lb]{\smash{{\SetFigFont{9}{10.8}{\familydefault}{\mddefault}{\updefault}{\color[rgb]{0,0,0}$x_2$}%
}}}}
\put(4145,-2942){\makebox(0,0)[lb]{\smash{{\SetFigFont{9}{10.8}{\familydefault}{\mddefault}{\updefault}{\color[rgb]{0,0,0}$x_1$}%
}}}}
\put(5168,-2814){\makebox(0,0)[lb]{\smash{{\SetFigFont{9}{10.8}{\familydefault}{\mddefault}{\updefault}{\color[rgb]{0,0,0}$x_3$}%
}}}}
\end{picture}%

%% file: detbspcurve.pstex_t
\begin{picture}(0,0)%
\includegraphics{detbspcurve.pstex}%
\end{picture}%
\setlength{\unitlength}{3947sp}%
\begingroup\makeatletter\ifx\SetFigFont\undefined%
\gdef\SetFigFont#1#2#3#4#5{%
  \reset@font\fontsize{#1}{#2pt}%
  \fontfamily{#3}\fontseries{#4}\fontshape{#5}%
  \selectfont}%
\fi\endgroup%
\begin{picture}(3919,1573)(4944,-4048)
\put(5294,-2695){\makebox(0,0)[lb]{\smash{{\SetFigFont{8}{9.6}{\familydefault}{\mddefault}{\updefault}{\color[rgb]{0,0,0}$x_1$}%
}}}}
\put(5586,-3569){\makebox(0,0)[lb]{\smash{{\SetFigFont{8}{9.6}{\familydefault}{\mddefault}{\updefault}{\color[rgb]{0,0,0}$l_5$}%
}}}}
\put(5410,-3103){\makebox(0,0)[lb]{\smash{{\SetFigFont{8}{9.6}{\familydefault}{\mddefault}{\updefault}{\color[rgb]{0,0,0}$l_1$}%
}}}}
\put(5819,-3220){\makebox(0,0)[lb]{\smash{{\SetFigFont{8}{9.6}{\familydefault}{\mddefault}{\updefault}{\color[rgb]{0,0,0}$l_2$}%
}}}}
\put(5877,-2578){\makebox(0,0)[lb]{\smash{{\SetFigFont{8}{9.6}{\familydefault}{\mddefault}{\updefault}{\color[rgb]{0,0,0}$x_3$}%
}}}}
\put(5935,-3745){\makebox(0,0)[lb]{\smash{{\SetFigFont{8}{9.6}{\familydefault}{\mddefault}{\updefault}{\color[rgb]{0,0,0}$x_2$}%
}}}}
\put(5235,-3336){\makebox(0,0)[lb]{\smash{{\SetFigFont{8}{9.6}{\familydefault}{\mddefault}{\updefault}{\color[rgb]{0,0,0}$l_4$}%
}}}}
\put(5586,-3394){\makebox(0,0)[lb]{\smash{{\SetFigFont{8}{9.6}{\familydefault}{\mddefault}{\updefault}{\color[rgb]{0,0,0}$l_3$}%
}}}}
\put(7743,-3161){\makebox(0,0)[lb]{\smash{{\SetFigFont{8}{9.6}{\familydefault}{\mddefault}{\updefault}{\color[rgb]{0,0,0}$h(x_1)$}%
}}}}
\put(4944,-3686){\makebox(0,0)[lb]{\smash{{\SetFigFont{8}{9.6}{\familydefault}{\mddefault}{\updefault}{\color[rgb]{0,0,0}$x_4$}%
}}}}
\end{picture}%

%% file: detbsp.pstex_t
\begin{picture}(0,0)%
\includegraphics{detbsp.pstex}%
\end{picture}%
\setlength{\unitlength}{3947sp}%
\begingroup\makeatletter\ifx\SetFigFont\undefined%
\gdef\SetFigFont#1#2#3#4#5{%
  \reset@font\fontsize{#1}{#2pt}%
  \fontfamily{#3}\fontseries{#4}\fontshape{#5}%
  \selectfont}%
\fi\endgroup%
\begin{picture}(5049,2094)(3826,-4123)
\put(4276,-2311){\makebox(0,0)[lb]{\smash{{\SetFigFont{11}{13.2}{\familydefault}{\mddefault}{\updefault}{\color[rgb]{0,0,0}$x_1$}%
}}}}
\put(8476,-3211){\makebox(0,0)[lb]{\smash{{\SetFigFont{11}{13.2}{\familydefault}{\mddefault}{\updefault}{\color[rgb]{0,0,0}$L_1$}%
}}}}
\put(5026,-2161){\makebox(0,0)[lb]{\smash{{\SetFigFont{11}{13.2}{\familydefault}{\mddefault}{\updefault}{\color[rgb]{0,0,0}$x_3$}%
}}}}
\put(3826,-3586){\makebox(0,0)[lb]{\smash{{\SetFigFont{11}{13.2}{\familydefault}{\mddefault}{\updefault}{\color[rgb]{0,0,0}$x_4$}%
}}}}
\put(5101,-3661){\makebox(0,0)[lb]{\smash{{\SetFigFont{11}{13.2}{\familydefault}{\mddefault}{\updefault}{\color[rgb]{0,0,0}$x_2$}%
}}}}
\put(7426,-2911){\makebox(0,0)[lb]{\smash{{\SetFigFont{11}{13.2}{\familydefault}{\mddefault}{\updefault}{\color[rgb]{0,0,0}$P_1=(0,0)$}%
}}}}
\put(7876,-3736){\makebox(0,0)[lb]{\smash{{\SetFigFont{11}{13.2}{\familydefault}{\mddefault}{\updefault}{\color[rgb]{0,0,0}$P_2=(1,-2)$}%
}}}}
\put(8176,-2611){\makebox(0,0)[lb]{\smash{{\SetFigFont{11}{13.2}{\familydefault}{\mddefault}{\updefault}{\color[rgb]{0,0,0}$(3,1)$}%
}}}}
\put(6526,-3361){\makebox(0,0)[lb]{\smash{{\SetFigFont{11}{13.2}{\familydefault}{\mddefault}{\updefault}{\color[rgb]{0,0,0}$(-1,-1.5)$}%
}}}}
\put(6976,-3811){\makebox(0,0)[lb]{\smash{{\SetFigFont{11}{13.2}{\familydefault}{\mddefault}{\updefault}{\color[rgb]{0,0,0}$L_2$}%
}}}}
\put(4426,-2836){\makebox(0,0)[lb]{\smash{{\SetFigFont{11}{13.2}{\familydefault}{\mddefault}{\updefault}{\color[rgb]{0,0,0}$l_1=2$}%
}}}}
\put(4951,-2986){\makebox(0,0)[lb]{\smash{{\SetFigFont{11}{13.2}{\familydefault}{\mddefault}{\updefault}{\color[rgb]{0,0,0}$l_2=\frac{1}{2}$}%
}}}}
\put(4651,-3211){\makebox(0,0)[lb]{\smash{{\SetFigFont{11}{13.2}{\familydefault}{\mddefault}{\updefault}{\color[rgb]{0,0,0}$l_3=1$}%
}}}}
\put(4651,-3436){\makebox(0,0)[lb]{\smash{{\SetFigFont{11}{13.2}{\familydefault}{\mddefault}{\updefault}{\color[rgb]{0,0,0}$l_5=1$}%
}}}}
\put(4126,-3136){\makebox(0,0)[lb]{\smash{{\SetFigFont{11}{13.2}{\familydefault}{\mddefault}{\updefault}{\color[rgb]{0,0,0}$l_4=1$}%
}}}}
\end{picture}%

%% file: largem04.pstex_t
\begin{picture}(0,0)%
\includegraphics{largem04.pstex}%
\end{picture}%
\setlength{\unitlength}{3947sp}%
\begingroup\makeatletter\ifx\SetFigFont\undefined%
\gdef\SetFigFont#1#2#3#4#5{%
  \reset@font\fontsize{#1}{#2pt}%
  \fontfamily{#3}\fontseries{#4}\fontshape{#5}%
  \selectfont}%
\fi\endgroup%
\begin{picture}(5424,2424)(2239,-3973)
\put(5026,-2836){\makebox(0,0)[lb]{\smash{{\SetFigFont{12}{14.4}{\familydefault}{\mddefault}{\updefault}{\color[rgb]{0,0,0}$l$}%
}}}}
\put(6826,-3886){\makebox(0,0)[lb]{\smash{{\SetFigFont{12}{14.4}{\familydefault}{\mddefault}{\updefault}{\color[rgb]{0,0,0}$S$}%
}}}}
\put(3076,-3811){\makebox(0,0)[lb]{\smash{{\SetFigFont{12}{14.4}{\familydefault}{\mddefault}{\updefault}{\color[rgb]{0,0,0}$C'$}%
}}}}
\put(5026,-3286){\makebox(0,0)[lb]{\smash{{\SetFigFont{12}{14.4}{\familydefault}{\mddefault}{\updefault}{\color[rgb]{0,0,0}$l'$}%
}}}}
\end{picture}%

%% file: stueck.pstex_t
\begin{picture}(0,0)%
\includegraphics{stueck.pstex}%
\end{picture}%
\setlength{\unitlength}{4144sp}%
\begingroup\makeatletter\ifx\SetFigFont\undefined%
\gdef\SetFigFont#1#2#3#4#5{%
  \reset@font\fontsize{#1}{#2pt}%
  \fontfamily{#3}\fontseries{#4}\fontshape{#5}%
  \selectfont}%
\fi\endgroup%
\begin{picture}(1620,741)(664,-544)
\end{picture}%

%% file: stringmov1.pstex_t
\begin{picture}(0,0)%
\includegraphics{stringmov1.pstex}%
\end{picture}%
\setlength{\unitlength}{3947sp}%
\begingroup\makeatletter\ifx\SetFigFont\undefined%
\gdef\SetFigFont#1#2#3#4#5{%
  \reset@font\fontsize{#1}{#2pt}%
  \fontfamily{#3}\fontseries{#4}\fontshape{#5}%
  \selectfont}%
\fi\endgroup%
\begin{picture}(4899,1787)(5239,-5201)
\put(8626,-4786){\makebox(0,0)[lb]{\smash{{\SetFigFont{10}{12.0}{\familydefault}{\mddefault}{\updefault}{\color[rgb]{0,0,0}$h(\Gamma')$}%
}}}}
\put(5626,-3886){\makebox(0,0)[lb]{\smash{{\SetFigFont{10}{12.0}{\familydefault}{\mddefault}{\updefault}{\color[rgb]{0,0,0}$h(e)$}%
}}}}
\put(5776,-5161){\makebox(0,0)[lb]{\smash{{\SetFigFont{10}{12.0}{\familydefault}{\mddefault}{\updefault}{\color[rgb]{0,0,0}(i)}%
}}}}
\put(6751,-4786){\makebox(0,0)[lb]{\smash{{\SetFigFont{10}{12.0}{\familydefault}{\mddefault}{\updefault}{\color[rgb]{0,0,0}$h(\Gamma')$}%
}}}}
\put(8176,-3961){\makebox(0,0)[lb]{\smash{{\SetFigFont{10}{12.0}{\familydefault}{\mddefault}{\updefault}{\color[rgb]{0,0,0}$h(e)$}%
}}}}
\put(8776,-5161){\makebox(0,0)[lb]{\smash{{\SetFigFont{10}{12.0}{\familydefault}{\mddefault}{\updefault}{\color[rgb]{0,0,0}(ii)}%
}}}}
\end{picture}%

%% file: contracted2.pstex_t
\begin{picture}(0,0)%
\includegraphics{contracted2.pstex}%
\end{picture}%
\setlength{\unitlength}{4144sp}%
\begingroup\makeatletter\ifx\SetFigFont\undefined%
\gdef\SetFigFont#1#2#3#4#5{%
  \reset@font\fontsize{#1}{#2pt}%
  \fontfamily{#3}\fontseries{#4}\fontshape{#5}%
  \selectfont}%
\fi\endgroup%
\begin{picture}(2816,2083)(2147,-1226)
\put(3016,-511){\makebox(0,0)[lb]{\smash{{\SetFigFont{10}{12.0}{\familydefault}{\mddefault}{\updefault}{\color[rgb]{0,0,0}$v_{1}$}%
}}}}
\put(2926,-61){\makebox(0,0)[rb]{\smash{{\SetFigFont{10}{12.0}{\familydefault}{\mddefault}{\updefault}{\color[rgb]{0,0,0}$v_{2}$}%
}}}}
\put(3061,119){\makebox(0,0)[lb]{\smash{{\SetFigFont{10}{12.0}{\familydefault}{\mddefault}{\updefault}{\color[rgb]{0,0,0}$v_{3}$}%
}}}}
\put(3466,389){\makebox(0,0)[lb]{\smash{{\SetFigFont{10}{12.0}{\familydefault}{\mddefault}{\updefault}{\color[rgb]{0,0,0}$v_{4}$}%
}}}}
\put(2881,-241){\makebox(0,0)[rb]{\smash{{\SetFigFont{10}{12.0}{\familydefault}{\mddefault}{\updefault}{\color[rgb]{0,0,0}$w_{1}$}%
}}}}
\put(2701,164){\makebox(0,0)[rb]{\smash{{\SetFigFont{10}{12.0}{\familydefault}{\mddefault}{\updefault}{\color[rgb]{0,0,0}$w_{2}$}%
}}}}
\put(3106,344){\makebox(0,0)[rb]{\smash{{\SetFigFont{10}{12.0}{\familydefault}{\mddefault}{\updefault}{\color[rgb]{0,0,0}$w_{3}$}%
}}}}
\put(3781,749){\makebox(0,0)[b]{\smash{{\SetFigFont{10}{12.0}{\familydefault}{\mddefault}{\updefault}{\color[rgb]{0,0,0}$\Gamma'$}%
}}}}
\put(2971,-1186){\makebox(0,0)[lb]{\smash{{\SetFigFont{10}{12.0}{\familydefault}{\mddefault}{\updefault}{\color[rgb]{0,0,0}(iii)}%
}}}}
\put(4366,-1186){\makebox(0,0)[lb]{\smash{{\SetFigFont{10}{12.0}{\familydefault}{\mddefault}{\updefault}{\color[rgb]{0,0,0}(iv)}%
}}}}
\end{picture}%

%% file: contractedf2.pstex_t
\begin{picture}(0,0)%
\includegraphics{contractedf2.pstex}%
\end{picture}%
\setlength{\unitlength}{4144sp}%
\begingroup\makeatletter\ifx\SetFigFont\undefined%
\gdef\SetFigFont#1#2#3#4#5{%
  \reset@font\fontsize{#1}{#2pt}%
  \fontfamily{#3}\fontseries{#4}\fontshape{#5}%
  \selectfont}%
\fi\endgroup%
\begin{picture}(5284,1191)(-51,-2845)
\put(1356,-1816){\makebox(0,0)[lb]{\smash{{\SetFigFont{10}{12.0}{\familydefault}{\mddefault}{\updefault}{\color[rgb]{0,0,0}$w_{1}$}%
}}}}
\put(1026,-2786){\makebox(0,0)[b]{\smash{{\SetFigFont{10}{12.0}{\familydefault}{\mddefault}{\updefault}{\color[rgb]{0,0,0}$v_{2}$}%
}}}}
\put(456,-2176){\makebox(0,0)[rb]{\smash{{\SetFigFont{10}{12.0}{\familydefault}{\mddefault}{\updefault}{\color[rgb]{0,0,0}$v_{1}$}%
}}}}
\put(2241,-2286){\makebox(0,0)[lb]{\smash{{\SetFigFont{10}{12.0}{\familydefault}{\mddefault}{\updefault}{\color[rgb]{0,0,0}$v_{2}$}%
}}}}
\put(1581,-2391){\makebox(0,0)[rb]{\smash{{\SetFigFont{10}{12.0}{\familydefault}{\mddefault}{\updefault}{\color[rgb]{0,0,0}$v_{1}$}%
}}}}
\put(1621,-2711){\makebox(0,0)[b]{\smash{{\SetFigFont{10}{12.0}{\familydefault}{\mddefault}{\updefault}{\color[rgb]{0,0,0}$w_{1}$}%
}}}}
\put(2776,-1811){\makebox(0,0)[b]{\smash{{\SetFigFont{10}{12.0}{\familydefault}{\mddefault}{\updefault}{\color[rgb]{0,0,0}$w_{1}$}%
}}}}
\put(3041,-2796){\makebox(0,0)[b]{\smash{{\SetFigFont{10}{12.0}{\familydefault}{\mddefault}{\updefault}{\color[rgb]{0,0,0}$v_{1}$}%
}}}}
\put(3371,-2231){\makebox(0,0)[lb]{\smash{{\SetFigFont{10}{12.0}{\familydefault}{\mddefault}{\updefault}{\color[rgb]{0,0,0}$v_{2}$}%
}}}}
\put(4726,-2246){\makebox(0,0)[lb]{\smash{{\SetFigFont{10}{12.0}{\familydefault}{\mddefault}{\updefault}{\color[rgb]{0,0,0}$v_{2}$}%
}}}}
\put(4406,-2516){\makebox(0,0)[lb]{\smash{{\SetFigFont{10}{12.0}{\familydefault}{\mddefault}{\updefault}{\color[rgb]{0,0,0}$v_{1}$}%
}}}}
\put(3806,-2641){\makebox(0,0)[b]{\smash{{\SetFigFont{10}{12.0}{\familydefault}{\mddefault}{\updefault}{\color[rgb]{0,0,0}$w_{1}$}%
}}}}
\put(1281,-1786){\makebox(0,0)[rb]{\smash{{\SetFigFont{10}{12.0}{\familydefault}{\mddefault}{\updefault}{\color[rgb]{0,0,0}$L_{1}$}%
}}}}
\put(1221,-2106){\makebox(0,0)[lb]{\smash{{\SetFigFont{10}{12.0}{\familydefault}{\mddefault}{\updefault}{\color[rgb]{0,0,0}$L_{2}$}%
}}}}
\put(2421,-1786){\makebox(0,0)[rb]{\smash{{\SetFigFont{10}{12.0}{\familydefault}{\mddefault}{\updefault}{\color[rgb]{0,0,0}$L_{1}$}%
}}}}
\put(2471,-2011){\makebox(0,0)[lb]{\smash{{\SetFigFont{10}{12.0}{\familydefault}{\mddefault}{\updefault}{\color[rgb]{0,0,0}$L_{2}$}%
}}}}
\put(3466,-2086){\makebox(0,0)[lb]{\smash{{\SetFigFont{10}{12.0}{\familydefault}{\mddefault}{\updefault}{\color[rgb]{0,0,0}$L_{2}$}%
}}}}
\put(3421,-1861){\makebox(0,0)[rb]{\smash{{\SetFigFont{10}{12.0}{\familydefault}{\mddefault}{\updefault}{\color[rgb]{0,0,0}$L_{1}$}%
}}}}
\put(4546,-1861){\makebox(0,0)[rb]{\smash{{\SetFigFont{10}{12.0}{\familydefault}{\mddefault}{\updefault}{\color[rgb]{0,0,0}$L_{1}$}%
}}}}
\put(4591,-2086){\makebox(0,0)[lb]{\smash{{\SetFigFont{10}{12.0}{\familydefault}{\mddefault}{\updefault}{\color[rgb]{0,0,0}$L_{2}$}%
}}}}
\end{picture}%

%% file: nocontractedA.pstex_t
\begin{picture}(0,0)%
\includegraphics{nocontractedA.pstex}%
\end{picture}%
\setlength{\unitlength}{4144sp}%
\begingroup\makeatletter\ifx\SetFigFont\undefined%
\gdef\SetFigFont#1#2#3#4#5{%
  \reset@font\fontsize{#1}{#2pt}%
  \fontfamily{#3}\fontseries{#4}\fontshape{#5}%
  \selectfont}%
\fi\endgroup%
\begin{picture}(733,1636)(860,-1460)
\put(891,-291){\makebox(0,0)[rb]{\smash{{\SetFigFont{10}{12.0}{\familydefault}{\mddefault}{\updefault}{\color[rgb]{0,0,0}$\check{w}_{2}=(0,1)$}%
}}}}
\put(1261,-511){\makebox(0,0)[lb]{\smash{{\SetFigFont{10}{12.0}{\familydefault}{\mddefault}{\updefault}{\color[rgb]{0,0,0}$\check{v}_{k}=\check{v}_{3}=(-1,2)$}%
}}}}
\put(1146,-741){\makebox(0,0)[lb]{\smash{{\SetFigFont{10}{12.0}{\familydefault}{\mddefault}{\updefault}{\color[rgb]{0,0,0}$\check{v}_{2}=(-1,1)$}%
}}}}
\put(1081,-1411){\makebox(0,0)[lb]{\smash{{\SetFigFont{10}{12.0}{\familydefault}{\mddefault}{\updefault}{\color[rgb]{0,0,0}$\check{v}_{1}=(-1,0)$}%
}}}}
\put(891,-971){\makebox(0,0)[rb]{\smash{{\SetFigFont{10}{12.0}{\familydefault}{\mddefault}{\updefault}{\color[rgb]{0,0,0}$\check{w}_{1}=(0,1)$}%
}}}}
\end{picture}%

%% file: stringS.pstex_t
\begin{picture}(0,0)%
\includegraphics{stringS.pstex}%
\end{picture}%
\setlength{\unitlength}{4144sp}%
\begingroup\makeatletter\ifx\SetFigFont\undefined%
\gdef\SetFigFont#1#2#3#4#5{%
  \reset@font\fontsize{#1}{#2pt}%
  \fontfamily{#3}\fontseries{#4}\fontshape{#5}%
  \selectfont}%
\fi\endgroup%
\begin{picture}(2131,699)(305,-523)
\put(721,-331){\makebox(0,0)[rb]{\smash{{\SetFigFont{10}{12.0}{\familydefault}{\mddefault}{\updefault}{\color[rgb]{0,0,0}$C_2$}%
}}}}
\put(991,-106){\makebox(0,0)[rb]{\smash{{\SetFigFont{10}{12.0}{\familydefault}{\mddefault}{\updefault}{\color[rgb]{0,0,0}$C_1$}%
}}}}
\put(1486, 29){\makebox(0,0)[b]{\smash{{\SetFigFont{10}{12.0}{\familydefault}{\mddefault}{\updefault}{\color[rgb]{0,0,0}$E'$}%
}}}}
\put(2206,-16){\makebox(0,0)[lb]{\smash{{\SetFigFont{10}{12.0}{\familydefault}{\mddefault}{\updefault}{\color[rgb]{0,0,0}$S$}%
}}}}
\put(1756,-286){\makebox(0,0)[lb]{\smash{{\SetFigFont{10}{12.0}{\familydefault}{\mddefault}{\updefault}{\color[rgb]{0,0,0}$E$}%
}}}}
\put(1261,-466){\makebox(0,0)[b]{\smash{{\SetFigFont{10}{12.0}{\familydefault}{\mddefault}{\updefault}{\color[rgb]{0,0,0}$E''$}%
}}}}
\end{picture}%